\documentclass[twocolumn]{autart}

\usepackage{graphics}
\usepackage{graphicx}
\usepackage{amssymb}
\usepackage{verbatim}
\usepackage{amsxtra}
\usepackage{epsfig}
\usepackage{float}
\usepackage[breaklinks = true]{hyperref} 
\usepackage{amsmath}
\usepackage{accents}			% This needs to be loaded after amsmath!!!!!!!!

\usepackage{mathtools}

\usepackage{latexsym}
\usepackage{ifthen}
\usepackage{color}
\usepackage{hyperref}
\usepackage{balance}
\usepackage{subeqnarray}
\usepackage{caption}
\usepackage{subcaption}
\usepackage[dvipsnames]{xcolor}

\usepackage{rgsEnvironments_auto}				%%
\usepackage{rgsMacros_auto}					%%

\usepackage{enumitem}
\usepackage{algorithm,algpseudocode}

\usepackage{ifthen}
\newboolean{Auto}
\setboolean{Auto}{false}  
\newcommand{\IfAuto}[2]{\ifthenelse{\boolean{Auto}}{{\color{black}#1}}{{\color{black}#2}}}

\newcommand{\OnlyForAuto}[1]{\IfAuto{#1}{}}
\newcommand{\NotForAuto}[1]{\IfAuto{}{#1}}

\usepackage{siunitx}

\usepackage{optidef}

\usepackage{pgfplots}

\usetikzlibrary{external}
\usepgfplotslibrary{external}
\usepackage{pgfplots}
\pgfplotsset{compat=newest,every axis/.append style={
                    label style={font=\footnotesize},
                    tick label style={font=\footnotesize}  
                    }}

\usepgfplotslibrary{fillbetween}

\DeclareMathOperator*{\interior}{int}
\DeclareMathOperator*{\closure}{cl}

\usepackage{amsmath}
\DeclareMathOperator*{\argmax}{arg\,max}

\newcommand{\Anorm}[1]{\left|#1\right|_{\mathcal{A}}}

\newcommand{\T}{\mathcal{T}}

\newcommand{\berk}[1]{{\color{blue}{#1}}}

\definecolor{mycolor}{RGB}{0,128,128}

\graphicspath{{./Figures/}}

\edef\endfrontmatter{%
  \unexpanded\expandafter{\endfrontmatter}% the current code
  \noexpand\endNoHyper % add \endNoHyper at the end to match \NoHyper
}

\begin{document}

%%%%%%%%%%%%%%%%%%%%%%%%%%%%%%%%%%%%%%%%%%%%%%%%%%%%%%%%%%%%%

\begin{frontmatter}
\runtitle{Semicontinuity Properties in Hybrid Dynamical Systems}  % Running title for regular 
                                              % papers but only if the title  
                                              % is over 5 words. Running title 
                                              % is not shown in output.

\title{Solutions and Reachable Sets of Hybrid Dynamical Systems: Semicontinuous Dependence on Initial Conditions, Time, and Perturbations\thanksref{footnoteinfo}} % Title, preferably not more 
                                                % than 10 words.

\thanks[footnoteinfo]{This paper was not presented at any IFAC 
meeting. Corresponding author B.~Alt{\i}n.}

\author{Berk Alt{\i}n}\ead{berkaltin@ucsc.edu},    % Add the  
\author{Ricardo G. Sanfelice}\ead{ricardo@ucsc.edu}  % (ead) as shown

%\address[berk]{Department of Electrical and Computer Engineering, University of California, Santa Cruz}
\address{Department of Electrical and Computer Engineering, University of California, Santa Cruz}

                             % keyword list or with the 
                                          % help of the Automatica 
                                          % keyword wizard

% Abstract of not more than 200 words.
\begin{abstract}
The sequential compactness afforded hybrid systems under mild regularity constraints guarantee outer/upper semicontinuous dependence of solutions on initial conditions and perturbations. For reachable sets of hybrid systems, this property leads to upper semicontinuous dependence with respect to initial conditions, time, and perturbations. Motivated by these results, we define a counterpart to sequential compactness and show that it leads to lower semicontinuous dependence of solutions on initial conditions and perturbations. In the sequel, it is shown that under appropriate assumptions, reachable sets of systems possessing this novel property depend lower semicontinuously on initial conditions, time, and perturbations. When those assumptions fail, continuous approximations of reachable sets turn out to be still possible. Necessary and sufficient conditions for the introduced property are given by a combination of geometric constraints, regularity assumptions, and tangentiality conditions. Further applications to simulations and optimal control of hybrid systems are discussed.
\end{abstract}

\end{frontmatter}

%%%%%%%%%%%%%%%%%%%%%%%%%%%%%%%%%%%%%%%%%%%%%%%%%%%%%%%%%%%%%

\section{Introduction}
\label{sec:intro}
The regularity of the map from initial conditions to solutions plays a key role in determining the structural properties of a dynamical system. As established in several textbooks on nonlinear continuous-time systems modeled as differential equations (e.g., \cite{khalil}), this map is continuous (i.e., solutions depend continuously on initial conditions) when the right-hand side of the differential equation is locally Lipschitz. Unfortunately, when the systems under consideration have constraints and/or set-valued dynamics, the problem becomes more challenging, and continuous dependence might be too much to ask. The challenges associated with constraints and set-valued dynamics are amplified for hybrid dynamical systems, since solutions from nearby initial conditions can have drastically different behavior due to the interaction between continuous and discrete dynamics. Nevertheless, in the last decade, it has been shown that semicontinuous dependence (in the upper sense) on initial conditions is enough to confer the set of solutions of a hybrid system with good enough structural properties to generate the following key results~\cite{hybridbook}: i) the basin of attraction of asymptotically stable compact sets are open; ii) asymptotic stability can be certified using a weak Lyapunov function via the invariance principle; and iii) asymptotic stability of a compact set is uniform.

Vaguely speaking, {\em upper semicontinuous dependence on initial conditions} is the property that given a nominal initial condition and a solution with perturbed initial condition, there exists a solution from the nominal initial condition that is close to the perturbed solution over a compact time horizon (in an appropriate sense). Generalizing this property to account for perturbations to the dynamics has enabled further results showcasing robustness of asymptotic stability and existence of smooth Lyapunov functions, thereby laying the foundations of a general modeling framework for hybrid systems (called \textit{hybrid inclusions}) that stands out on its emphasis on a robust stability theory~\cite{hybridsurvey}. For hybrid inclusions, upper semicontinuous dependence is guaranteed when the set of solutions is sequentially compact; that is, the limit of a convergent sequence of solutions is also a solution. Importantly, as shown originally in~\cite{GOEBEL2006573}, sequential compactness is implied when the data defining the system certifies mild regularity assumptions called the \textit{hybrid basic conditions}; see~\cite[Ch.~6]{hybridbook}. Among the many properties afforded hybrid systems satisfying the hybrid basic conditions is robustness to a large class of perturbations, which has been shown to include singular perturbations~\cite{sanfelice2011} and temporal perturbations such as delays~\cite{delay}. Importantly, the robustness properties inherited from the hybrid basic conditions has been put to use in diverse application areas such as robot locomotion~\cite{5535144}, power systems~\cite{theunisse}, event-triggered control and vehicle platooning~\cite{DolPlo_TITS17a}, network control~\cite{HeiPos_AUT18a}, and desynchronization~\cite{7870654}, among many others. Despite these very promising developments, conditions guaranteeing continuous dependence on initial conditions and perturbations have remained elusive, with the notable exceptions of~\cite{cai}, where a hybrid Filippov-Wa\.{z}ewski relaxation theorem is developed, and~\cite{BROUCKE2002149,semantics,980123,1582902}, where continuous dependence has been studied for hybrid automata without state perturbations, with the results often assuming the guards/domains to be differentiable manifolds and requiring open sets (we do not assume the sets defining a hybrid inclusion to be smooth). In this paper, we address this issue, with the following specific contributions.
\begin{enumerate}[label={\arabic*)},leftmargin=*]
	\item We provide necessary and sufficient conditions for the solutions of hybrid systems to depend {\em lower semicontinuously on initial conditions and perturbations}; i.e., given a nominal solution and a perturbation of the initial condition and the dynamics, there exists a solution of the perturbed system from the perturbed initial condition close to the nominal solution. 
	\item The necessary and sufficient conditions requiring partial knowledge of the solutions, we show how results from viability theory, which depend solely on the data of the system, can be utilized to check that the necessary and sufficient conditions hold.
	\item For hybrid systems whose solutions depend upper and lower semicontinuously on initial conditions and perturbations, we show that their reachable sets depend continuously on initial conditions, perturbations, and time under appropriate assumptions, and how they can be continuously approximated even when those assumptions fail.
\end{enumerate}

The technical developments here are made possible through a natural counterpart of the aforementioned sequential compactness property, introduced in an earlier version of this work~\cite{cdc2020}. In plain words, the introduced property guarantees the existence of a sequence of solutions convergent to any given nominal solution. In Section~\ref{sec:solutionset} of this article, we expand on this property by generalizing it to account for perturbations, and show that solutions of systems possessing this property depend lower semicontinuously on initial conditions and perturbations. To guarantee the existence of a sequence of (perturbed) solutions convergent to a given nominal one, sufficient (and necessary) conditions are given in Section~\ref{sec:nomiwp}, which have not appeared in~\cite{cdc2020}. As mentioned, these require partial knowledge of solutions, in particular, those of the underlying continuous-time system. However, the viability arguments we put forth in Section~\ref{sec:viab}, which depend only on the data of the hybrid system, ensure that the sufficient conditions hold.

As noted above, an important application of the developed theory concerns reachable sets of hybrid systems, in particular, their continuous dependence on initial conditions, perturbations, and time. For reachable sets, although there has been significant interest on the development of computational algorithms (e.g.~\cite{7526475,ALTHOFF2010233,geret}), continuity properties have not attracted much attention in the literature, even for differential \IfAuto{inclusions (see the technical report accompanying this article~\cite{arxiv}).}{inclusions.\footnote{To the best of our knowledge, attention has been restricted mostly to showing upper semicontinuity with respect to the initial conditions, e.g. Theorem 1 in Section 2 of~\cite[Ch.~2]{aubincellina} or~\cite[Proposition~3.5.5]{viability}. Corollary~5.3.3 in \cite{viability} and Corollary~10.4.2 in~\cite{aubin} establish Lipschitz continuity in the case of Lipschitz inclusions whose solutions do not reach the boundary of the domain of the flow map. We extend these to the more general hybrid case and show continuity with respect to initial conditions, time and perturbations. Lower semicontinuity results have been stated for differential inclusions \textit{without constraints} in~\cite{eps,donch}.}} We show in Section~\ref{sec:reachprop} that while there are unique challenges associated with the hybrid setting and even for differential inclusions with constraints, upper/lower semicontinuous dependance of solutions on initial conditions and perturbations can be used to continuously approximate reachable sets. Unlike the initial results in~\cite{cdc2020}, the continuous approximation results here do not assume any knowledge of solutions, while allowing state perturbations.

\NotForAuto{
The significance of our results are underlined by their connections to a) continuous approximations of solutions and reachable sets via discretization, and b) regularity properties of value functions and optimal controls, especially those arising in model predictive control~\cite{acc}. Existing hybrid systems discretization/simulation frameworks guarantee that simulated solutions are not too far off from true solutions by virtue of upper semicontinuous dependence on initial conditions and state perturbations, thereby practically preserving asymptotic stability. However, they do not guarantee that simulated solutions faithfully reproduce \textit{all solutions}, and as such, cannot reliably verify safety/invariance properties and accurately compute optimal controls. The results of this work set the foundation for these missing links, and will enable a consistent theory of discretization for hybrid systems guaranteeing that true solutions can be recovered as the discretization step size is decreased, in the spirit of~\cite{wolenski}, where the finite-time reachable set of a Lipschitz differential inclusion is given by an exponential formula corresponding to the forward Euler method; see also~\cite{dontchev} for a more general discussion.\footnote{Another related result by Frankowska~\cite[Theorem~3.5.6]{viability}, which has appeared previously as~\cite[Corollary~2.6]{frankowska}, establishes the right-hand side of the differential inclusion as the infinitesimal generator of the reachable set mapping.} One particular avenue where such a theory would be particularly useful is model predictive control~\cite{acc}, where quantifiable connections between abstract formulations and their computationally tractable counterparts~\cite{9029729} can be made.}

We highlight that the results derived in this article encompass constrained differential inclusions as a special case, for which results are limited and have mostly considered solutions that live in an open set~\cite{aubincellina,viability}. Moreover, due to our inherent assumption of set-valued dynamics, the developed tools are applicable to systems with inputs recast as autonomous hybrid inclusions, by virtue of Filippov's lemma (e.g.~\cite[Corollary~23.4]{clarke}). With respect to~\cite{cai}, the results here have a generally different flavor. Aside from our consideration of perturbations, the metric we use to quantify closeness between solutions is consistent with~\cite{hybridbook}, and the viability conditions we employ are more general than those in~\cite{cai}, with the exception of convexity of the flow map (as that work pertains to a relaxation theorem, although some of the viability conditions therein also assume convexity); e.g. compare~\ref{item:R4} in Theorem~\ref{thrm:viab} with condition (4) in \cite[Assumption~4.3]{cai}. The sufficient (and necessary) conditions in Section~\ref{sec:nomiwp} can also be interpreted to be more general then the relaxation conditions in~\cite{cai}; for example, the semicontinuity conditions on the discrete dynamics are less demanding.

Our work has some similarities to that of~\cite{semantics} on semi-computability of solutions on hybrid automata, which relies on notions of upper and lower semi-continuity of hybrid automata. Aside from the choice of modeling framework, some of the differences with this work can be summarized as follows: a) we parametrize reachable sets by initial condition, ordinary time, \textit{and the number of jumps} (and at times, perturbations); b) though motivated by a theory of computation, we focus on (semi)continuity properties of solutions and reachable sets, which at times lead to stronger results than those achievable for computability; for instance, in contrast with the results of~\cite{semantics} that show that solutions are generally not upper and lower semi-computable at the same time, we prove in \IfAuto{\cite[Proposition 8]{arxiv}}{Proposition~\ref{prop:ionwp}} that a hybrid system can be both \textit{nominally inner and outer well-posed}, notions similar to lower/upper semicontinuity in~\cite{semantics}; c) the notions of \textit{well-posedness} of hybrid systems in our paper, which resemble the notions of semicontinuity of hybrid automata in~\cite{semantics}, are not stated directly in terms of regularity conditions; rather, we provide regularity conditions that are sufficient for such well posedness. We would also like to note that while some of the challenges associated with discontinuities discussed in~\cite{semantics} have also been identified in our work, we provide a concrete way of continuously approximating reachable sets even at the points of discontinuities using class-$\mathcal{K}$ estimates.

\section{Preliminaries}
\label{sec:background}
We use~$\reals$ to denote real numbers,~$\realsgeq$ nonnegative reals, and~$\naturals$ nonnegative integers. The 2-norm is denoted~$|.|$. For a given pair of sets~$S_1,S_2$,~$S_1\subset S_2$ indicates~$S_1$ is a subset of~$S_2$, not necessarily proper. Let~$\A\subset\reals^n$ be nonempty. The distance of a vector~$x\in\reals^n$ to the set~$\A$ is~${\Anorm{x}:=\inf_{a\in\A}|x-a|}$. The closed unit ball in~$\reals^n$ centered at the origin is denoted~$\ball$, and~$\A+r\ball$ is the set of all~$x$ such that~$|x-a|\leq r$ for some~$a\in\A$. The closure, interior, and boundary of a set~${S\subset\reals^n}$ are denoted~${\closure S}$,~${\interior S}$, and~${\partial S}$. The domain of a set-valued mapping~$M:S\rightrightarrows \reals^m$, denoted~$\dom M$, is the set of all~$x\in S$ such that~$M(x)$ is nonempty. Given a set~$S'\subset S$,~$M|_{S'}$ denotes the restriction of~$M$ to~$S'$. A continuous function~$\alpha:\realsgeq\to\realsgeq$ is a class-$\mathcal{ K}$ function if it is strictly increasing and~$\alpha(0)=0$.

\subsection{Hybrid Inclusions and Hybrid Arcs}

We consider hybrid systems in the form
\begin{equation}
	\HS
	\left\{
	\begin{aligned}
		\dot{x}		&\in F(x)		& x&\in C\\
		x^+				&\in G(x)		& x&\in D,
	\end{aligned}
	\right.
	\label{eq:H}
\end{equation}
where the \textit{flow map}~${F:\reals^n\rightrightarrows \reals^n}$ defines the continuous-time evolution~(\textit{flows}) of the state~$x\in\reals^n$ on the \textit{flow set}~${C\subset\dom F}$, and the \textit{jump map}~${G:\reals^n\rightrightarrows \reals^n}$ defines the discrete transitions~(\textit{jumps}) of~$x$ on the \textit{jump set}~$D\subset\dom G$ \cite{hybridbook}. To refer to the hybrid system in~\eqref{eq:H} and define its data~$(C,F,D,G)$, at times, we use the notation~$\HS=(C,F,D,G)$. A special case of~\eqref{eq:H} is when~$D$ is empty, which corresponds to a constrained continuous-time system, denoted compactly as~$(C,F)$.

Solutions of the hybrid system~$\HS$ belong to a class of functions called \textit{hybrid arcs} and are parametrized by hybrid time~$(t,j)$, where~$t\in\realsgeq$ denotes the ordinary time and~$j\in\nats$ denotes the number of jumps. A function~$x$ mapping a subset of~$\realsgeq\times\naturals$ to~$\reals^n$ is a hybrid arc if 1) its domain, denoted~$\dom x$, is a \textit{hybrid time domain}, and 2) it is locally absolutely continuous on each connected component of~$\dom x$. Formally, a set $E \subset\realsgeq\times\naturals$ is a hybrid time domain if for every~$(T,J)\in E$, there exists a nondecreasing sequence~$\{t_j\}_{j=0}^{J+1}$ with~$t_0=0$ such that~$E \cap \left([0,T]\times \{0,1,\dots,J\}\right) = \cup_{j=0}^{J}\left([t_j,t_{j+1}]\times \{j\}\right)$, which implies~$t_{j+1}=T$.  Equivalently, $E$ is a (possibly finite) disjoint union of intervals~$[s_0,s_1],[s_1,s_2],\dots$, where $0=s_0\leq s_1\leq\dots$, with the last interval (if it exists) possibly open to the right and/or unbounded. A hybrid arc~$x$ with~$x(0,0)\in \closure(C)\cup D$ is a solution of the hybrid system~$\HS$ if the following hold:
\begin{itemize}
	\item for every $j\geq 0$ with nonempty~$I^j:=\{t: (t,j)\in\dom x\}$, $x(t,j)\in C$ for all~$t\in \interior I^j$ and~$\dot{x}(t,j)\in F(x(t,j))$ for almost all~$t\in I^j$;
	\item for all~$(t,j)\in\dom $ such that~$(t,j+1)\in\dom x$, $x(t,j)\in D$ and~$x(t,j+1)\in G(x(t,j))$.
\end{itemize}
In the case of a continuous-time system~$(C,F)$, for simplicity, we omit the jump index~$j$ and parametrize the solutions only by ordinary time~$t$.

A hybrid arc~$x$ is called 1) \textit{trivial} if its domain is a singleton and \textit{nontrivial} otherwise; 2) \textit{complete} if its domain is unbounded; and 3) \textit{continuous} if it is nontrivial and its domain is connected (i.e., $\dom x\subset \realsgeq\times\{0\}$). It is called \textit{bounded} if its range is bounded. It is said to \textit{escape to infinity at hybrid time~$(T,J)$} (or have \textit{finite escape time}) if~$x(t,J)$ tends to infinity as~$t$ tends to~$T$ from the left. If the domain of~$x$ is compact, we say that~$(T,J)\in\dom x$ is the \textit{terminal (hybrid) time} of~$x$ if~$t\leq T$ and~$j\leq J$ for all~$(t,j)\in\dom x$. Similarly,~$T$ is referred to as the \textit{terminal ordinary time} of~$x$. A hybrid arc~$x$ with terminal time~$(T,J)$ is said to \textit{terminate on a set~$S$} if~$x(T,J)\in S$. The same terminology is used for hybrid arcs that are solutions of the hybrid system~$\HS$; e.g., a solution~$x$ of~$\HS$ is bounded if its range is bounded.

A solution~$x$ of the hybrid system~$\HS$ is \textit{maximal} if it cannot be extended to another solution. The notation~$\sol_{\HS}(S)$ refers to the set of all maximal solutions~$x$ of~$\HS$ originating from~$S$ (i.e.~${x(0,0)\in S}$ for every~$x\in\sol_{\HS}(S)$), and~$\sol_{\HS}:=\sol_{\HS}(\reals^n)$. If every~$x\in\sol_{\HS}(S)$ is bounded or complete, we say that~$\HS$ is \textit{pre-forward complete from~$S$}. We say that~$t$ is a \textit{jump time} of~$x$ if there exists~$j$ such that~$(t,j),(t,j+1)\in\dom x$. The (nontrivial) solution~$x$ is said to \textit{begin with a flow} if~$[0,\varepsilon]\times\{0\}\subset \dom x$ for some~$\varepsilon>0$. It is said to \textit{begin with a jump} if~$(0,1)\in\dom x$. For the hybrid system~$\HS=(C,F,D,G)$, flows are said to be possible at~$x_0$ if there exists a solution~$x$ of~$\HS$ originating from~$x_0$ that begins with a flow. The set of points where flows are possible is the set
\begin{equation}
	\widetilde{C}:=\{x_0: \exists x\in\sol_{\HS}(x_0), \varepsilon>0 \text{ s.t. } (\varepsilon,0)\in\dom x\}.
\label{eq:setflow}
\end{equation}

We use~$x$ interchangeably to denote both the state of~$\HS$ (or points in~$\reals^n$) and solutions of~$\HS$, and explicitly specify what it refers to when necessary.

\subsection{Set-Valued Analysis Background}

Let~$S\subset\reals^n$,~$x\in \closure S$, and consider a set-valued mapping~$M:S\rightrightarrows \reals^m$. The \textit{inner limit} of~$M$ as~$x'$ tends to~$x$, denoted~$\liminf_{x'\to x}M(x')$, is the set of all~$y$ such that for any sequence~$\{x_i\}_{i=0}^{\infty}\in S$ convergent to~$x$, there exist~$\imath\geq 0$ and a sequence~$\{y_i\}_{i=\imath}^{\infty}$ convergent to~$y$ such that~$y_i\in M(x_i)$ for all~$i\geq \imath$. The \textit{outer limit} of~$M$ as~$x'$ tends to~$x$, denoted~$\limsup_{{x'\to x}}M(x')$, is the set of all~$y$ for which there exists a sequence~$\{x_i\}_{i=0}^{\infty}\in S$ convergent to~$x$ and a sequence~$\{y_i\}_{i=0}^{\infty}$ convergent to~$y$ such that~$y_i\in M(x_i)$ for all~$i\geq 0$. When the inner and outer limits (as~$x'$ tends to~$x$) are equal, the limit of~$M$ as~$x'$ tends to~$x$, denoted~$\lim_{x'\to x}M(x')$, is defined to be equal to them. Limits of sequences of sets are defined in the same manner. Let~$X\subset S$ and~$x\in \closure X$. Then, the mapping~$M$ is said to be \textit{inner semicontinuous} (respectively, \textit{outer semicontinuous}) at~$x$ relative to~$X$ if the inner (respectively, outer) limit of~$M|_{X}$ as~$x'$ tends to~$x$ contains (respectively, is contained in) $M(x)$.\NotForAuto{\footnote{Alternatively, in a more general context, these properties can be stated using  hyperspace topologies~\cite{beer}.}} It is said to be \textit{continuous} at~$x$ relative to~$X$ if it is both inner and outer semicontinuous at~$x$ relative to~$X$. In addition,~$M$ is \textit{locally bounded at~$x\in X$ relative to~$X$} if there exists~$\varepsilon>0$ such that the set~$M((x+\varepsilon\ball)\cap X)$ is bounded. If these properties hold for all~$x\in X$, we drop the qualifier ``at~$x$'', and if~$X=S$, we drop the qualifier ``relative to~$X$''. Also, the mapping~$M$ is \textit{Lipschitz} on~$X$ if it has nonempty values on~$X$ and there exists~$L\geq 0$ such that~$M(x)\subset M(x')+L|x-x'|\ball$ for every~$x,x'\in X$. \IfAuto{The definitions here follow their counterparts in~\cite{rockafellarwets}.}{

The definitions of set convergence, semicontinuity, and local boundedness here follow~\cite[Definitions~4.1,~5.4, and~5.14]{rockafellarwets}. For locally bounded set-valued maps with closed values, outer semicontinuity is also equivalent to the property commonly known as \textit{upper semicontinuity}~\cite[Definition~1.4.1]{aubin}, see~\cite[Lemma~5.15]{hybridbook}. Inner semicontinuity coincides with the property commonly known as \textit{lower semicontinuity}~\cite[Definition~1.4.2]{aubin}.}

%%%%%%%%%%%%%%%%%%%%%%%%%%%%%%%%%%%%%%%%%%%%%%%%%%%%%%%%%%%%%%%%%%%%%%%%%%%%%%

\subsection{Graphical Convergence and Outer Well-Posedness}

Let~$\{x_i\}_{i=0}^{\infty}$ be a sequence of hybrid arcs. The sequence~$\{x_i\}_{i=0}^{\infty}$ is \textit{locally eventually bounded} if for any~$\tau\geq 0$, there exist~$\imath\geq 0$ and a compact set~$K$ such that~$x_i(t,j)\in K$ for every~${i\geq \imath}$ and~$(t,j)\in\dom x_i$ with~$t+j\leq\tau$. It is said to \textit{converge graphically} to a mapping~$M:\realsgeq\times\nats\rightrightarrows\reals^n$ if the set limit of the sequence~$\{\gph x_i	\}_{i=0}^{\infty}$ equals~$\gph M$, where~$\gph $ denotes the graph of a set-valued mapping. The mapping~$M$ is the \textit{graphical limit} of~$\{x_i\}_{i=0}^{\infty}$, and has closed graph since the set limit is always closed. See~\cite[Chapter~5]{hybridbook} for details. A related concept that allows one to quantify closeness in the hybrid setting is called~$(\tau,\varepsilon)$-closeness~\cite[Definition~5.23]{hybridbook}.

\begin{defn}
Given~$\tau\geq 0$ and $\varepsilon>0$, two hybrid arcs~$x$ and~$x'$ are said to be~$(\tau,\varepsilon)$-close if
	\begin{itemize}[leftmargin=*]
		\item for every~$(t,j)\in\dom x$ satisfying~$t+j\leq\tau$, there exists~$(t',j)\in\dom x'$ such that $|t-t'|<\varepsilon$ and $|x(t,j)-x'(t',j)|<\varepsilon$;
		\item for every~$(t',j')\in\dom x'$ satisfying~$t'+j'\leq\tau$, there exists~$(t,j')\in\dom x$ such that $|t'-t|<\varepsilon$ and $|x'(t',j')-x(t,j')|<\varepsilon$.
	\end{itemize}
\end{defn}

The graphical limit~$x$ of a locally eventually bounded and graphically convergent sequence of hybrid arcs~$\{x_i\}_{i=0}^{\infty}$ is defined on a hybrid time domain. In fact,~$\dom x$ is precisely the set limit of the corresponding sequence of hybrid time domains. If~$x$ is single valued on~$\dom x$, then for each~$j\geq 0$, the value of~$x$ at the beginning and end of the interval~$I^{j}:=\{t:(t,j)\in \dom x\}$ is given by the limit of corresponding points of the sequence~$\{x_i\}_{i=0}^{\infty}$, \IfAuto{see~\cite[Lemma 2]{arxiv}.}{as shown below. The proof follows easily from~\cite[Theorem~5.25]{hybridbook} and is not included.

\begin{lem}[Convergence of Jumps]
\label{lem:term}
Let~$\{x_i\}_{i=0}^{\infty}$ be a locally eventually bounded sequence of hybrid arcs with closed domains. Suppose that the sequence graphically converges to a hybrid arc~$x$, and let~$t_j^i:=\min I^j_i$ and~$t_{j+1}^i:=\sup I^j_i$ for every~$i,j\geq 0$ such that the interval~$I^j_i:=\{t:\exists (t,j)\in\dom x_i\}$ is nonempty. Given~$j\geq 0$, let~$I^j:=\{t:\exists (t,j)\in\dom x\}$. Then, the following hold.
\begin{itemize}[leftmargin=*]
	\item The interval~$I^j$ is nonempty if and only if there exists~$\imath\geq 0$ such that~$t_j^i$ is defined for all~$i\geq \imath$ and the sequence~$\{t_j^i\}_{i=\imath}^{\infty}$ is bounded. Moreover, if there exists such~$\imath$, the sequence~$\{(t_j^i,x_i(t_j^i,j))\}_{i=\imath}^{\infty}$ converges to $(t_j,x(t_j,j))$, where~$t_j:=\min I^j$.
	\item The interval~$I^j$ is nonempty and bounded if and only if there exists~$\imath\geq 0$ such that~$t_{j+1}^i$ is defined for all~$i\geq\imath$ and the sequence~$\{t_{j+1}^i\}_{i=\imath}^{\infty}$ is bounded. Moreover, if there exists such~$\imath$, the sequence $\{(t_{j+1}^i,x_i(t_j^i,j))\}_{i=\imath}^{\infty}$ converges to $(t_{j+1},x(t_{j+1},j))$, where~$t_{j+1}:=\max I^j$.
\end{itemize}
\end{lem}
}

In~\cite{hybridbook}, \textit{nominally well-posed} hybrid systems are defined to have a graphical convergence property that can be interpreted as outer semicontinuous dependence of solutions on initial conditions: for a nominally well-posed system~$\HS$, the graphical limit~$x$ of a locally eventually bounded graphically convergent sequence~$\{x_i\}_{i=0}^{\infty}$ of solutions is itself a solution, with~$x(0,0)=\lim_{i\to\infty}x_i(0,0)$. For this reason, given a set~$S$, in this article, hybrid systems possessing the property outlined in~\cite[Definition~6.2]{hybridbook} for sequences of solutions whose initial conditions converge to a point in~$S$ are said to be nominally \textit{outer} well-posed on~$S$. The definition is recalled in Appendix~\ref{sec:defs} for completeness. Hybrid systems possessing a similar property in the presence of vanishing state perturbations (c.f.~\cite[Definition~6.29]{hybridbook}) are said to be \textit{outer well-posed} on~$S$;\footnote{For the various notions of well-posedness in the paper, for simplicity, we omit the qualifier ``on $S$'' when~$S=\reals^n$. Also, we say ``at $x_0$'' instead of ``on $S$'' if $S=\{x_0\}$ for some~$x_0$. Note this and the various notions throughout the paper apply to continuous-time systems when the jump set is empty and to discrete-time systems when the flow set is empty.} see Definition~\ref{def:owp} in Appendix~\ref{sec:defs}, which relies on Definition~\ref{def:rho} below. The change in terminology is to accommodate the counterpart of these definitions guaranteeing inner semicontinuous dependence of solutions on initial conditions (and perturbations).

By definition, every outer well-posed hybrid system is nominally outer well-posed. Outer well-posedness, being a property motivated by robustness, admits a fairly general class of perturbations called~$\rho$-perturbations~(\cite[Definition 6.27]{hybridbook}), recalled below.

\begin{defn}[$\rho$-Perturbation]
\label{def:rho}
Given a hybrid system~$\HS=(C,F,D,G)$ and a function~$\rho:\reals^n\to\realsgeq$, the~$\rho$-perturbation of~$\HS$ is the hybrid system~$\HS^{\rho}$ with data~$(C^{\rho},F^{\rho},D^{\rho},G^{\rho})$, where $C^{\rho}=\{x: (x+\rho(x)\ball)\cap C\neq \varnothing\}$, $D^{\rho}=\{x: (x+\rho(x)\ball)\cap D\neq \varnothing\}$, and 
	\begin{equation*}
		\begin{aligned}
			F^{\rho}(x)&=\closure(\con F((x+\rho(x)\ball)\cap C))+\rho(x)\ball,\\
			G^{\rho}(x)&=\{z: z\in y+\rho(y)\ball, y\in G((x+\rho(x)\ball)\cap D)\},
		\end{aligned}
	\end{equation*}
for all~$x\in\reals^n$, where $\con$ denotes the convex hull. Moreover, given any~$\delta\in(0,1)$,~$\HS^{\delta\rho}$ denotes the~$\tilde{\rho}$-perturbation of~$\HS$, where~$\tilde{\rho}$ is the function~$x\mapsto\delta\rho(x)$.
\end{defn}

Observe that given~$\delta_1\leq \delta_2$, every solution of $\HS^{\delta_1\rho}$ is a solution of~$\HS^{\delta_2\rho}$. Outer well-posedness (and hence, nominal outer well-posedness) can be verified by checking that the data of the system satisfies the \textit{hybrid basic conditions}; \cite[Assumption~6.5 and Theorem~6.8]{hybridbook}. 

\begin{thm}
\label{thrm:nomout}
A hybrid system~$\HS=(C,F,D,G)$ is outer well-posed if the following hold.
\begin{enumerate}[label={(A\arabic*)},leftmargin=*]
	\item \label{item:A1}	The sets~$C$ and~$D$ are closed.
	\item \label{item:A2}	The flow map~$F$ is locally bounded and outer semicontinuous relative to~$C$, and~$C\subset \dom F$. Furthermore, for every~${x\in C}$, the set~${F(x)}$ is convex.
	\item \label{item:A3}	The jump map~$G$ is locally bounded and outer semicontinuous relative to~$D$, and~$D\subset \dom G$.
\end{enumerate}
\end{thm}

\begin{exmp}[Bouncing Ball]
\label{ex:bouncingball}
Consider a ball bouncing vertically on a horizontal flat surface. When modeled as a point-mass with height~$x_1$ and velocity~$x_2$, the motion of the ball can be represented by the hybrid system $\HS$ in~\eqref{eq:H} with state~$x:=(x_1,x_2)$, where $C=\{x:x_1\geq 0\}$, $D=\{x:x_1=0, x_2\leq 0\}$, and for every~$x\in\reals^2$, $F(x)=(x_2,-\gamma)$ and $G(x)=(0,-\lambda x_2)$. Here,~$\gamma>0$ is the gravitational acceleration and~$\lambda\in[0,1]$ is the coefficient of restitution. It is straightforward to see that this system satisfies the hybrid basic conditions and is therefore outer well-posed, since in the case of single-valued maps, \ref{item:A2}-\ref{item:A3} are equivalent to continuity of the flow map~$F$ on the flow set~$C$ and the jump map~$G$ on the jump set~$D$.
\end{exmp}

\section{Inner Well-Posedness of Hybrid Systems}
\label{sec:solutionset}
We now introduce ``inner'' counterparts of the outer well-posedness concepts in Section~\ref{sec:solutionset}.

\subsection{Nominally Inner Well-Posed Hybrid Systems}

The notion of nominal outer well-posedness comprises two mutually exclusive cases: when the sequence of solutions in question is locally bounded, the graphical limit~$x$ is a solution with closed graph (since the set limit is always closed) that is necessarily bounded or complete, and therefore has closed domain. If it is not locally bounded, then the graphical limit leads to a (maximal) solution~$x$ that escapes to infinity. See~\cite{hybridbook} for details. Consequently, nominally inner well-posed hybrid systems, which can be interpreted as those hybrid systems whose solutions depend \textit{inner} semicontinously on initial conditions, are defined as follows.

\begin{defn}[Nominal Inner Well-Posedness]
\label{def:mydef}
A hybrid system~$\HS=(C,F,D,G)$ is said to be nominally inner well-posed on a set~$S$ if for every solution~$x$ of $\HS$ originating from~$S$, the following holds.
\begin{enumerate}[label={($\star$)},leftmargin=*]
	\item \label{item:star} Given any sequence~$\{\xi_{i}\}_{i=0}^{\infty}\in\closure(C)\cup D$ convergent to~$x(0,0)$, for every~$i\geq 0$, there exists a solution~$x_i$ of~$\HS$ originating from~$\xi_i$ such that
	\begin{enumerate}[label={(\alph*)},leftmargin=*]
		\item \label{item:stara} if~$x$ is complete or bounded with~$\dom x$ closed, then the sequence of solutions~$\{x_i\}_{i=0}^{\infty}$ is locally eventually bounded and graphically convergent to~$x$;
		\item \label{item:starb} if~$x$ escapes to infinity at hybrid time~$(T,J)$, then the sequence of solutions~$\{x_i\}_{i=0}^{\infty}$ is not locally eventually bounded but graphically convergent to a mapping~$M$ such that~$x=M|_{\dom M\cap([0,T)\times\{0,1,\dots,J\})}$.
	\end{enumerate}
\end{enumerate} 
\end{defn}

\begin{exmp}[Thermostat]
\label{ex:therm}
Consider the hybrid system model~$\HS=(C,F,D,G)$ of a closed-loop thermostat with state~$x:=(z,q)$, where~$z\in\reals$ is the room temperature and~$q\in\{0,1\}$ is a binary variable denoting whether thermostat is on or off. Given desired minimum and maximum temperatures~$z_{\min}$ and~$z_{\max}$, respectively, suppose that~$z_o<z_{\min}<z_{\max}<z_o+z_{\Delta}$, where~$z_o$ is the natural temperature of the room and~$z_{\Delta}>0$ is the capacity of the heater to raise the temperature. The data of the system is given as~$C=\{x:z\geq z_{\min}, q=0 \text{ or } z\leq z_{\max},q=1\}$, $D=\{z\leq z_{\min},q=0 \text{ or }z\geq z_{\max},q=1\}$, and~$F(x)=(-z+z_o+qz_{\Delta},0)$ and~$G(x)=(z,1-q)$ for all~$x\in\reals^2$. Note that maximal solutions of this system are unique and with some abuse of notation, let~$x=(z,q)$ be the maximal solution of~$\HS$ originating from~$(z_{\min},0)$. Consider the sequence in~\ref{item:star}, and let~$x_i$ be the maximal solution from~$\xi_i$ for each~$i\geq 0$. For each~$i\geq 0$, if~$\xi_i=(z_{\min}+\varepsilon_i,0)$ for some~$\varepsilon_i>0$, $x_i$ flows until reaching~$x(0,0)=(z_{\min},0)$ and then jumps to~$(z_{\min},1)$, otherwise, it jumps immediately and flows to reach~$(z_{\min},1)$. Using this fact and uniqueness of maximal solutions, and also observing that~$\lim_{i\to \infty}T_i=0$, where~$T_i$ satisfies~$x_i(T_i,1)=(z_{\min},1)$ for all~$i\geq 0$, it can be concluded that~$\{x_i\}_{i=0}^{\infty}$ is locally eventually bounded and graphically convergent to~$x$. 
\end{exmp}

Necessary and sufficient conditions for nominal inner well-posedness that depend only on the data of the system are provided via the results in Sections~\ref{sec:nomiwp} and~\ref{sec:viab}. \IfAuto{For a pre-forward complete nominally inner and outer well-posed system, solutions depend upper/lower semicontinuously on initial conditions, uniformly over compact sets of initial conditions~\cite[Proposition~3.13]{cdc2020} (see also~\cite[Section~3]{arxiv} for a more nuanced discussion). In general, nominal outer and inner well-posedness do not imply each other. However, when maximal solutions are unique and complete, they turn out to be equivalent \cite[Proposition 8]{arxiv}.}{Under an absolute continuity assumption during flows,\footnote{The results in~\cite{cdc2020} assume that~$x$ has a closed graph. Fortunately, if a maximal solution~$x$ is absolutely continuous at each interval of flow, then its graph is necessarily closed. Otherwise, one may encounter the pathological case of~$x$ exhibiting infinitely fast changes. For example, the \textit{locally} absolutely continuous function~$x(t)=\sin(1/(1-t))$ for all~$t\in[0,1)$ is a maximal solution of~$\dot{x}\in(-\infty,\infty)$. Nominal outer well-posedness excludes such behavior.} the graphical convergence property in~\ref{item:star} is equivalent to an alternative formulation using~$(\tau,\varepsilon)$-closeness for bounded or complete solutions~\cite[Propositions~3.8 and~3.9]{cdc2020}. Moreover, for a hybrid system~$\HS$, lower semicontinuous dependence on initial conditions\footnote{See (\cite[Definition~3.10]{cdc2020}, and~\cite{cai} for a similar definition with a slightly different terminology and measure of closeness.} implies nominal inner well-posedness (\cite[Theorem~3.11]{cdc2020}). The reverse implication is true if~$\HS$ is nominally outer well-posed and pre-forward complete (\cite[Theorem~3.12]{cdc2020}). Hence, solutions of a pre-forward complete nominally inner and outer well-posed system depend upper/lower semicontinuously on initial conditions, uniformly over compact sets of initial conditions~\cite[Proposition~3.13]{cdc2020}.

In general, nominal outer and inner well-posedness do not imply each other. However, when maximal solutions are unique and complete, they turn out to be equivalent.

\begin{prop}
\label{prop:ionwp}
Let~$\HS$ be a hybrid system. Given an initial condition~$x_0$, suppose that~$\HS$ has a unique maximal solution originating from~$x_0$, and there exists~$r> 0$ such that every maximal solution of~$\HS$ originating from~$x_0+r\ball$ is complete. Then,~$\HS$ is nominally inner well-posed at~$x_0$ if it is nominally outer well-posed at~$x_0$. The reverse implication is true if~$\HS$ has a unique maximal solution originating from~$x_0'$ for every~$x'_0\in x_0+r\ball$.
\end{prop}
\begin{pf}
For the first implication, let~$x$ be the unique maximal solution from~$x_0$. Take~$\varepsilon>0$ and~$\tau\geq\varepsilon+1$. Since~$x$ is complete, there exists~$\delta\in(0,r)$ such that the following holds: for any maximal solution~$x'$ with~$x'(0,0)\in x_0+\delta\ball$, there exists a truncation of~$x$, say~$\tilde{x}$, such that~$\tilde{x}$ and~$x'$ are~$(2\tau,\varepsilon)$-close by \cite[Proposition~3.13]{cdc2020}. Take any maximal solution~$x'$ originating from~$x_0+r\ball$, which is complete, and let~$\tilde{x}$ be the corresponding truncation of~$x$ such that~$\tilde{x}$ and~$x'$ are~$(2\tau,\varepsilon)$-close. Let $(T',J'):=\argmax \{t+j: (t,j)\in\dom x', t+j\leq 2\tau\}$ and note that~$T'+J'>2\tau-1$. Then, by closeness between~$\tilde{x}$ and~$x'$, there exists~$T>T'-\varepsilon$ such that~$(T,J')\in\dom \tilde{x}$. Hence,~$T+J'>2\tau-\varepsilon-1$. Given~$(t',j)\in\dom x'$ with~$t'+j\leq \tau$, there exists~$t$ such that~$(t,j)\in\dom \tilde{x}\subset\dom x$, $|t'-t|<\varepsilon$, and~$|x'(t',j)-\tilde{x}(t,j)|=|x'(t',j)-x(t,j)|<\varepsilon$. Now pick~$(t,j)\in\dom x$ such that~$t+j\leq \tau$. Since~$\tau\geq \varepsilon+1$, $t+j \leq T+J'$ and therefore~$(t,j)\in\dom\tilde{x}$. Thus, there exists~$t'$ such that~$(t',j)\in\dom x'$, $|t-t'|<\varepsilon$, and~$|\tilde{x}(t,j)-x'(t',j)|=|x(t,j)-x'(t',j)|<\varepsilon$. Consequently,~$x$ and~$x'$ are~$(\tau,\varepsilon)$-close. This implies that for every~$x'_0\in (x_0+\delta\ball)$, there exists a maximal solution~$x'$ originating from~$x'_0$ such that~$x$ and~$x'$ are~$(\tau,\varepsilon)$-close, so by \cite[Proposition~3.8]{cdc2020},~\ref{item:star} holds. By considering truncations of~$x$,~$\HS$ is nominally inner well-posed at~$x_0$.

For the reverse implication, let~$x$ be the unique maximal solution from~$x_0$, which is complete by assumption, and has closed domain. Let~$\{x_i\}_{i=1}^{\infty}$ be a graphically convergent sequence of solutions of~$\HS$ with~$\lim_{i\to\infty}x_i(0,0)=x(0,0)$. By nominal inner well-posedness,~\ref{item:star} holds, so using \cite[Proposition~3.9]{cdc2020}, it can be concluded in a similar manner that for every~$\varepsilon>0$ and~$\tau\geq 0$, there exists~$\delta>0$ such that the following holds: for every maximal solution~$x'$ originating from~$x_0+\delta\ball$,~$x$ and~$x'$ are~$(\tau,\varepsilon)$-close, where. Then, by~\cite[Theorem~5.25]{hybridbook} and due to uniqueness of solutions, the graphical limit of the sequence~$\{x_i\}_{i=0}^{\infty}$ is precisely~$x$.
\end{pf}

} In the absence of uniqueness, similar conclusions cannot be reached. A simple counterexample is the nominally outer well-posed differential inclusion~$\dot{x}=F(x)$ on~$\reals$, where~$F(x)=1$ if~$x>0$, $F(x)=-1$ if~$x<0$, and~$F(0)=[-1,1]$, for which the solution~$x(t)=0 $ for all~$t\geq 0$ fails~\ref{item:star} in Definition~\ref{def:mydef}. The need for completeness is demonstrated by the following example.

\begin{exmp}[Planar Continuous-Time System]
\label{ex:planar}
Consider the continuous-time system~$(C,F)$ with flow set~$C=\{x=(x_1,x_2): x_1x_2=0, x_1\geq 0\}$ and flow map given as $F(x)=(1,0)$ for all~$x\in\reals^2$. Nominal outer well posedness of this system can be verified directly by Theorem~\ref{thrm:nomout}. However, this system is not nominally inner well-posed at the origin since maximal solutions from outside the~$x_1$-axis are all trivial, while the continuous maximal solution from the origin is complete.
\end{exmp}

\subsection{Inner Well-Posed Perturbations of Hybrid Systems}

Recall that outer well-posedness of a hybrid system generalizes nominal outer well-posedness so that the effects of state perturbations can be scrutinized. Inspired by this, given a hybrid system~$\HS$, roughly speaking, we refer to a family of parametrized hybrid systems $\{\HS_{\delta}=(C_{\delta},F_{\delta},F_{\delta},G_{\delta})\}_{\delta\in(0,1)}$ as an \textit{inner well-posed perturbation of~$\HS$} if an analogue of the graphical convergence property in \ref{item:star} of Definition~\ref{def:mydef} holds for sequences of solutions of the parametrized family. For a more compact notation, when referring to such families of systems, we omit the subscript~$\delta\in(0,1)$.

\begin{defn}[Inner Well-Posed Perturbations]
\label{def:iwppert}
A family of hybrid systems~$\{\HS_{\delta}=(C_{\delta},F_{\delta},D_{\delta},G_{\delta})\}$ is said to be an inner well-posed perturbation of a hybrid system~$\HS$ on a set~$S$ if~$S\cap(\closure(C)\cup D)\subset\liminf_{\delta\to 0}\closure(C_{\delta})\cup D_{\delta}$, and for every solution~$x$ of $\HS$ originating from~$S$, the following holds:
\begin{enumerate}[label={($\diamond$)},leftmargin=*]
	\item \label{item:asterisk} given any sequence~$\{\delta_i\}_{i=0}^{\infty}\in(0,1)$ convergent to zero and any sequence~$\{\xi_{i}\}_{i=0}^{\infty}$ convergent to~$x(0,0)$ with~$\xi_i\in\closure(C_{\delta_i})\cup D_{\delta_i}$ for all~$i\geq 0$, for every~$i\geq 0$, there exists a solution~$x_i$ of~$\HS_{\delta_i}$ originating from~$\xi_i$ such that~\ref{item:stara} and \ref{item:starb} in Definition~\ref{def:mydef} hold.
\end{enumerate}
\end{defn}

Note that because~$S\cap(\closure(C)\cup D)\subset\liminf_{\delta\to 0}\closure(C_{\delta})\cup D_{\delta}$, \ref{item:asterisk} is not a vacuous statement. That is, given the sequence~$\{\delta_i\}_{i=0}^{\infty}\in(0,1)$, there exists a sequence~$\{\xi_{i}\}_{i=0}^{\infty}$ convergent to~$x(0,0)$ with~$\xi_i\in\closure(C_{\delta_i})\cup D_{\delta_i}$ for all~$i\geq 0$.

When~$\HS$ is nominally inner well-posed, trivially, the family of hybrid systems $\{\HS_{\delta}\}$ satisfying $\HS_{\delta}=\HS$ for all~$\delta>0$ is an inner well-posed perturbation of~$\HS$. The primary motivation for inner well-posed perturbations is recovering inner semicontinuous dependence of solutions (in the sense of Definition~\ref{def:mydef}) for systems that are not nominally inner well-posed.
\begin{exmp}[Planar System Revisited]
Given the continuous-time system~$(C,F)$ in Example~\ref{ex:planar}, consider the family of continuous-time systems~$\{(C_{\delta},F_{\delta})\}$, where for every~$\delta\in(0,1)$,~$C_{\delta}:=\{x: x_1\geq 0, |x_2|\leq \delta\}$ and~$F_{\delta}(x):=(1,-\delta x_2)$ for all $x\in\reals^2$. Since~$C_{\delta}\supset C$ for all~$\delta\in(0,1)$, it is clear that $\closure C\subset \liminf_{\delta\to 0} \closure C_{\delta}$. Moreover, for any solution of~$(C,F)$ originating outside the origin, the graphical convergence property in~\ref{item:asterisk} holds trivially since the only nontrivial maximal solution of~$(C,F)$ originates from the origin. Now, observing that given any~$\delta>0$ and any initial condition~$\xi\in C_{\delta}$, the unique maximal solution of~$(C_{\delta},F_{\delta})$ from~$\xi=(\xi_1,\xi_2)$ is given as~$x'(t)=(\xi_1+t,\xi_2\exp(-\delta t))$ for all~$t\geq 0$, it can be shown that the convergence property in~\ref{item:asterisk} holds for any solution of~$(C,F)$ from the origin, so~$\{(C_{\delta},F_{\delta})\}$ is an inner well-posed perturbation of~$(C,F)$. 
\end{exmp}

\NotForAuto{
The next two results generalize \cite[Propositions 3.8 and~3.9]{cdc2020}. The proofs are similar and and are omitted for brevity.
\begin{prop}
\label{prop:iscgraphical0pert}
Let $x$ be a solution of a hybrid system~$\HS=(C,F,D,G)$ with closed graph. Given a family of hybrid systems~$\{\HS_{\delta}=(C_{\delta},F_{\delta},D_{\delta},G_{\delta})\}$, suppose that for every~$\varepsilon>0$ and~$\tau\geq 0$, there exist~$\bar{\delta},r>0$ such that the following holds: for any~$\delta\in(0,\bar{\delta}]$ and~$x'_0\in (x(0,0)+r\ball)\cap (\closure(C_{\delta})\cup D_{\delta})$, there exists a solution~$x'$ of~$\HS_{\delta}$ originating from~$x'_0$ such that~$x$ and~$x'$ are~$(\tau,\varepsilon)$-close. Then,~\ref{item:asterisk} holds.
\end{prop}
%\begin{proof}
%Consider the sequence $\{(k,1/k)\}_{k=1}^{\infty}$ and let $\{\bar{\delta}_k\}_{k=1}^{\infty}$ be a strictly decreasing positive sequence such that for every~$k\geq 1$ the following holds: for all~$\delta\leq \bar{\delta}=:r$ and~$x'_0\in (x(0,0)+\bar{\delta}_k\ball)\cap (\closure(C_{\delta})\cup D_{\delta})$, there exists a solution~$x'$ of~$\HS_{\delta}$ originating from~$x'_0$ such that $x$ and $x'$ are~$(k,1/k)$-close. Let $\{i'_k\}_{k=1}^{\infty}$ be a strictly increasing positive sequence such that given any~$k\geq 1$, $|\xi_i-x_0|\leq \bar{\delta}_k$ and~$\delta_i\leq\bar{\delta}_k$ for all $i\geq i'_k$. For every~$k\geq 1$ and every~$i\in\{i'_k,i'_k+1,\dots,i'_{k+1}-1\}$, pick a solution $x'_i$ of $\HS_{\delta_i}$ originating from $\xi_i$ such that $x'_i$ and $x$ are $(k,1/k)$-close. By~\cite[Theorem~5.25]{hybridbook}, the sequence $\{x'_i\}_{i=1}^{\infty}$ is graphically convergent to $x$. In the case of a bounded or complete~$x$, local eventual boundedness of $\{x'_i\}_{i=1}^{\infty}$ follows from the fact that for any~$\tau\geq 0$, the set of all~$x(t,j)$ with~$t+j\leq\tau$ is bounded. {\color{red}In addition, convergence of~$\{\dom x_i\}_{i=0}^{\infty}$ to~$\dom x$ and~$\{\length \dom x_i\}_{i=0}^{\infty}$ to~$\length \dom x$ follow via~\cite[Example~5.19]{hybridbook} and~\cite[Theorem~5.25]{hybridbook}.} Otherwise, if~$x$ escapes to infinity at hybrid time~$(T,J)$, it is obvious that $\{x'_i\}_{i=1}^{\infty}$ is not locally eventually bounded.
%\end{proof}
\begin{prop}
\label{prop:graphicalisc0pert}
Let~$x$ be a bounded or complete solution of a hybrid system~$\HS=(C,F,D,G)$ with closed domain, and suppose that given a family of hybrid systems~$\{\HS_{\delta}=(C_{\delta},F_{\delta},D_{\delta},G_{\delta})\}$, \ref{item:asterisk} holds. Then, for every~$\varepsilon>0$ and~$\tau\geq 0$, there exist~$\bar{\delta},r>0$ such that the following holds: for any~$\delta\in(0,\bar{\delta}]$ and~$x'_0\in (x(0,0)+r\ball)\cap (\closure(C_{\delta})\cup D_{\delta})$, there exists a solution~$x'$ of~$\HS_{\delta}$ originating from~$x'_0$ such that~$x$ and~$x'$ are~$(\tau,\varepsilon)$-close.
\end{prop}
%\begin{proof}
%If the conclusion of the proposition were false, there would exist a sequence~$\{\delta_i\}_{i=0}^{\infty}$ and~$\{\xi_i\}_{i=1}^{\infty}\in \closure(C_{\delta_i})\cup D_{\delta_i}$ such that for every~$i\geq 1$,~$\delta_i\leq 1/i$,~$\xi_i\in x(0,0)+(1/i)\ball$, and no solution~$x_i$ of~$\HS_{\delta_i}$ satisfying~$x_i(0,0)=\xi_i$ is such that~$x$ and~$x_i$ are~$(\tau,\varepsilon)$-close. Since~$x$ is bounded or complete and the domain of~$x$ is closed, the graph of~$x$ is closed. For each~$i\geq 1$, pick a solution~$x_i$ of~$\HS_{\delta_u}$ with initial condition~$x_i(0,0)=\xi_i$ such that the sequence~$\{x_i\}_{i=1}^{\infty}$ is locally eventually bounded and graphically convergent to~$x$. Then, by~\cite[Theorem~5.25]{hybridbook}, there exists~$\imath\in\nats$ such that for every~$i\geq \imath$,~$x$ and~$x_i$ are~$(\tau,\varepsilon)$-close, which is a contradiction.
%\end{proof}
}

At times, we work with family of hybrid systems that can be ``upper bounded'' by the~$\rho$-perturbation (recall Definition~\ref{def:rho}) of a hybrid system, formalized below.
\begin{defn}[Domination by a $\rho$-Perturbation]
A family of hybrid systems~$\{\HS_{\delta}=(C_{\delta},F_{\delta},D_{\delta},G_{\delta})\}$ is said to be dominated by the $\rho$-perturbation of a hybrid system~$\HS$ if for every~$\delta\in(0,1)$, $C_{\delta}\subset C^{\delta\rho}$, $F_{\delta}(x)\subset F^{\delta\rho}(x)$ for all~$x\in \reals^n$, $D_{\delta}\subset D^{\delta\rho}$, and $G_{\delta}(x)\subset G^{\delta\rho}(x)$ for all~$x\in \reals^n$, where~$(C^{\delta\rho},F^{\delta\rho},D^{\delta\rho},G^{\delta\rho})$ is the data of the $\delta\rho$-perturbation of~$\HS$; see Definition~\ref{def:rho}.
\end{defn}

\IfAuto
{
Similar to nominally inner-well posed systems, if an inner well-posed perturbation $\{\HS_{\delta}\}$ of a pre-forward complete outer well-posed~$\HS$ is dominated by a $\rho$-perturbation of $\HS$ for some continuous~$\rho$, then solutions of $\HS$ depend lower semicontinuously on initial conditions and perturbations~\cite[Proposition 15]{arxiv}.
}
{
If an inner well-posed perturbation~$\{\HS_{\delta}\}$ of a pre-forward complete~$\HS$ is dominated by a~$\rho$-perturbation of~$\HS$ for some continuous~$\rho$, then solutions of~$\HS$ depend lower semicontinuously on initial conditions and perturbations, generalizing~\cite[Proposition~3.13]{cdc2020}. The proof is very similar and thus not included.

\begin{prop}
\label{prop:contout}
Let~$\HS$ be a hybrid system, and given a compact set~$K$, suppose that~$\HS$ is outer well-posed on~$K$ and pre-forward complete from~$K$. Let~$\{\HS_{\delta}=(C_{\delta},F_{\delta},D_{\delta},G_{\delta})\}$ be an inner well-posed perturbation of~$\HS$ on~$K$, and suppose that~$\{\HS_{\delta}\}$ is dominated by a~$\rho$-perturbation of~$\HS$ for some continuous function~$\rho$. Then, for all~$\varepsilon>0$ and~$\tau\geq 0$, there exist~$r,\bar{\delta}>0$ such that the following holds: for every~$x\in\sol_{\HS}(K)$, $\delta\in(0,\bar{\delta}]$, and~$x'_0\in(x(0,0)+r\ball)\cap (\closure(C_{\delta})\cup D_{\delta})$ there exists a solution~$x'$ of~$\HS_{\delta}$ originating from~$x'_0$ such that~$x$ and~$x'$ are~$(\tau,\varepsilon)$-close.
\end{prop}

\begin{remark}
\label{rem:notvac}
For Proposition~\ref{prop:contout} to be meaningful, it should be possible to choose scalars~$r>0$ and~$\bar{\delta}>0$ such that the set~$(x(0,0)+r\ball)\cap (\closure(C_{\delta})\cup D_{\delta})$ is nonempty for all~$\delta\in(0,\bar{\delta}]$. This is indeed the case, which can be observed using the alternative definition of inner semicontinuity in~\cite[Proposition~5.12]{rockafellarwets}.
\end{remark}
} 

\section{In-Depth Look at Inner Well-Posedness}
\label{sec:nomiwp}
Let~$\HS$ be a nominally inner well-posed system and suppose that~$x$ is a solution of~$\HS$ that begins with a (nontrivial) flow and then jumps when it reaches the jump set. Directly by definition and \IfAuto{\cite[Lemma 2]{arxiv}}{recalling Lemma~\ref{lem:term}} and~\cite[Proposition~3.9]{cdc2020}, there must be a neighborhood~$S$ of~$x(0,0)$ such that flows are possible everywhere on the set~$S\cap\closure (C)$. Moreover, solutions originating from~$S\cap\closure (C)$ that begin with a flow must be able to reach the jump set approximately in the same amount of ordinary time as~$x$. This section formalizes such statements and derives necessary and sufficient conditions for nominal inner well-posedness and inner well-posed perturbations.

\subsection{Nominally Inner Well-Posed Systems}

The requirement that~$\HS$ must have solutions that reach the jump set approximately in the same amount of ordinary time as~$x$, postulated above, can be formalized by generalizing nominal inner well-posedness to account for \textit{terminal constraints}. As such, given a target set~$X$, we look at solutions of~$\HS$ that terminate on~$X$. {Recalling the definition of ``termination'' in Section~\ref{sec:background}, such a solution $x$ is not necessarily  maximal but has compact domain, and satisfy~$x(T,J)\in X$, where~$(T,J)$ is the terminal time of~$x$.} This generalization bears some similarities to the notions of ``relaxation with constraints and a target''~\cite{cai} and ``viability with a target'' \cite{quinveliov,aubinimpulse}, the latter of which is motivated by optimal control problems~\cite{quincampoix}. 

\begin{defn}
\label{def:mydef2}
A hybrid system~$\HS=(C,F,D,G)$ is said to be nominally inner well-posed on a set~$S$ with terminal constraint~$X$ if for every solution~$x$ of $\HS$ originating from~$S$ and terminating on~$X$, the following holds:
\begin{enumerate}[label={($\star\star$)},leftmargin=*]
	\item \label{item:star2} given any sequence~$\{\xi_{i}\}_{i=0}^{\infty}\in\closure(C)\cup D$ convergent to~$x(0,0)$, there exists~$\imath\geq 0$ such that for every~$i\geq \berk{\imath}$, there exists a solution~$x_i$ of~$\HS$ originating from~$\xi_i$ and terminating on~$X$ such that the sequence~$\{x_i\}_{i=\berk{\imath}}^{\infty}$ is locally eventually bounded and graphically convergent to~$x$.
\end{enumerate} 
\end{defn}

Graphical convergence and local boundedness of the sequence~$\{x_i\}_{i=\berk{\imath}}^{\infty}$ in~\ref{item:star2} implies that the sequence of terminal times (respectively, points) of~$\{x_i\}_{i=0}^{\infty}$ converges to the terminal time (respectively, point) of~$x$; \IfAuto{\cite[Lemma 2]{arxiv}}{see Lemma~\ref{lem:term}}. Consequently, given a hybrid system that is nominally inner well-posed with terminal constraint~$X$, for any solution~$x$ that terminates on~$X$, there exist solutions from a neighborhood of~$x(0,0)$ terminating on~$X$ that are close to~$x$ in the sense of~$(\tau,\varepsilon)$-closeness, as in~\cite[Proposition~3.9]{cdc2020}. 

Trivially, a nominally inner well-posed hybrid system~$\HS=(C,F,D,G)$ is nominally inner well-posed with terminal constraint~$\reals^n$ or $\closure(C)\cup D\cup G(D)$. Note that \ref{item:starb} in Definition~\ref{def:mydef2} is irrelevant for this implication since solutions escaping to infinity in finite time do not ``terminate'' according to our definition of ``termination''. Of course, more interesting examples abound. \IfAuto{For example, the closed-loop thermostat model~$\HS$ in Example~\ref{ex:therm} can be observed to be nominally inner well-posed with terminal constraint $X=[z_{\min},z_{\max}]\times\{0,1\}$, directly by definition of solutions. }{

\begin{exmp}[Thermostat Revisited]
Consider the closed-loop thermostat model~$\HS$ in Example~\ref{ex:therm}. By Theorem~\ref{thrm:nomout} and Proposition~\ref{prop:ionwp}, this system is nominally inner well-posed. In addition, directly by definition of solutions, it can be observed that every solution converges to the set~$X=[z_{\min},z_{\max}]\times\{0,1\}$ in finite time. Due to uniqueness of maximal solutions, this fact can be used to show nominal inner well-posedness with terminal constraint~$ X$.
\end{exmp}}

\NotForAuto{
\begin{remark}
In general, nominal inner well-posedness and nominal inner well-posedness with terminal constraints do not imply each other: the continuous-time system $(C,F)$ with~$C=\reals$ and~$F(x)=0$ for all~$x\in C$ is nominally inner well-posed but not nominally inner well-posed with terminal constraint $X=\{0\}$. The discontinuous differential equation $\dot{x}=F(x)$ on~$C=\mathbb{R}$ with~$F(x)=-x$ on~$X=(-1,1)$ and~$F(x)=x$ on~$C$ is not nominally inner well-posed at~$\{-1,1\}$ but is nominally inner well-posed with terminal constraint $X$.
\end{remark}}

\OnlyForAuto{In general, nominal inner well-posedness and nominal inner well-posedness with terminal constraints do not imply each other (see~\cite[Remark 4.1]{arxiv}). }Sufficient conditions for nominal inner well-posedness require nominal inner-well posedness (with the jump set as the terminal constraint) for the flows of~$\HS$, along with restrictions on the geometry of the flow and jump sets, and inner semicontinuity of the jump map. 

\begin{thm}[Sufficient Conditions]
\label{thrm:nomiwp}
Given a hybrid system~$\HS=(C,F,D,G)$, let~$\widetilde{C}$ be the set of all points where flows are possible, defined in~\eqref{eq:setflow}. Then,~$\HS$ is nominally inner well-posed if the following hold:
\begin{enumerate}[label={(B\arabic*)},leftmargin=*]
	\item \label{item:N1} for every~$x\in\widetilde{C}$,~$(x+r\ball)\cap(\closure(C)\cup D)\subset\widetilde{C}$ for some~$r>0$;
	\item \label{item:N2} for every~$x\in D$,~$(x+r\ball)\cap(\closure(C)\cup D)\subset\widetilde{C} \cup D$ for some~$r>0$;
	\item \label{item:N3} the continuous-time system $(C,F)$ is nominally inner well-posed;
	\item \label{item:N4} the continuous-time system $(C,F)$ is nominally inner well-posed with terminal constraint~$D$;
	\item \label{item:N5} the jump map~$G$ is inner semicontinuous relative to~$D$;
	\item \label{item:N6} the mapping~$\widetilde{G}:\reals^n\rightrightarrows\reals^n$, defined below, is inner semicontinuous relative to~$D$:
	\begin{equation}
		\widetilde{G}(x):=G(x)\cap(\widetilde{C}\cup D) \quad \forall x\in \reals^n.
	\label{eq:restrict}
	\end{equation}
\end{enumerate}
\end{thm}

The next example demonstrates the application of the sufficient conditions. The main idea behind the conditions of Theorem~\ref{thrm:nomiwp} is that inner semicontinuous approximation of flows is possible only by flows, while jumps can be approximated by jumps or a ``short'' flows followed by jumps\IfAuto{, as shown in the upcoming Example~\ref{ex:great}.}{. In particular, given a bounded solution~$x$ that begins with a (nontrivial) flow, Lemma~\ref{lem:term} implies that inner semicontinuous approximation (in the sense of Definition~\ref{def:mydef}) is possible only if given a sequence~$\{\xi_{i}\}_{i=0}^{\infty}\in\closure(C)\cup D$ convergent to~$x(0,0)$, there exists~$\imath\geq 0$ such that for every~$i\geq \imath$, there exists a continuous solution~$x_i$ originating from~$\xi_i$. This conclusion leads to~\ref{item:N1} and \ref{item:N3}. Condition \ref{item:N4} further ensures that if~$x$ has a jump at ordinary time~$T$, the approximating sequence~$\{x_i\}_{i=0}^{\infty}$ can be selected such that~$x_i$ has a jump for large enough~$i$. On the other hand, since jumps can be approximated by flows (see Example~\ref{ex:bbnowp} below), \ref{item:N2} only requires that either a jump or a flow is possible from nearby any initial condition in the jump set. In the case of a solution~$x$ that begins with a jump, \ref{item:N5} ensures that~$x$ can be inner semicontinuously approximated by solutions that also begin with a jump. Finally, if~$x$ begins with a jump and then either flows or jumps, \ref{item:N6} ensures that~$x$ can be inner semicontinuously approximated by solutions that begin with a jump and then either flows or jumps\footnote{For example, the unique maximal solution of the bouncing ball from the origin, which exhibits only jumps, is approximated in the sense of Definition~\ref{def:mydef} by any sequence of maximal solutions whose initial conditions tend to the origin. Provided the initial conditions are nonzero, maximal solutions of such sequences exhibit flows.}}

\begin{exmp}[Bouncing Ball Revisited]
\label{ex:bbnowp}
Consider the bouncing ball model in Example~\ref{ex:bouncingball} and an initial condition~$x=(x_1,x_2)$. For this system, flows are possible if and only if~$x_1>0$, or $x_1=0$ and~$x_2> 0$. Hence,~\ref{item:N1} holds. Similarly, since jumps are possible if and only if~$x_1=0$ and~$x_2\leq 0$, \ref{item:N2} holds. Moreover, one can show that \ref{item:N3}-\ref{item:N4} hold by using the closed-form solutions; see~\cite[Example 2.12.]{hybridbook}. Inner semicontinuity of the jump map~$G$ relative to the jump set~$D$ is obvious since the jump map is continuous, so \ref{item:N5} holds. For \ref{item:N6}, it suffices to note that the mapping~$\widetilde{G}$ in \eqref{eq:restrict} is equal to~$G$ on~$D$, as~$\widetilde{C}\cup D= C$ and~$G(D)\subset C$. Hence, this system is nominally inner well-posed.
\end{exmp}

The proof of Theorem~\ref{thrm:nomiwp} is postponed until the end of the section, as it is a special case of the upcoming Theorem~\ref{thrm:nomiwppert}. For the time being, we note that while some of the sufficient conditions outlined for nominal inner well-posedness are also necessary, the same cannot be said of others. In particular, inner semicontinuity of the jump map and the restricted jump map need not be true at points where flows are possible.

\begin{thm}[Necessary Conditions]
\label{prop:nec}
Given a nominally inner well-posed hybrid system~$\HS$ with data~$(C,F,D,G)$, let~$\widetilde{C}$ be the set of all points where flows are possible. Then,~\ref{item:N1}-\ref{item:N4} and the following hold:
\begin{enumerate}[label={(B\arabic*')},leftmargin=*]
\setcounter{enumi}{4}
\item \label{item:N5'} for every~$x\in D\cap \closure(D\backslash\widetilde{C})$, the jump map~$G$ is inner semicontinuous at~$x$ relative to~$D\backslash\widetilde{C}$.
\item \label{item:N6'} for every~$x\in D\cap \closure(D\backslash\widetilde{C})$, the mapping~$\widetilde{G}$ in~\eqref{eq:restrict} is inner semicontinuous at~$x$ relative to~$D\backslash\widetilde{C}$.
\end{enumerate}
\end{thm}
\NotForAuto{
\begin{pf}
Let~$x$ be a solution of~$\HS$ such that~\ref{item:star} holds. Necessity of~\ref{item:N1} and~\ref{item:N3} follow by assuming~$x$ to be continuous and noting that the sequence~$\{x_{i}\}_{i=0}^{\infty}$ in~\ref{item:star} would then have to be continuous for large~$i$ by Lemma~\ref{lem:term}. Similar arguments show necessity of~\ref{item:N4} by letting~$x$ to be a solution that begins with a flow, jumps once, and then terminates. Now let~$x$ be a solution that begins with a jump. Then,~\ref{item:N2} follows directly, and assuming the sequence~$\{\xi_{i}\}_{i=0}^{\infty}$ in~\ref{item:star} is in~$D\backslash\widetilde{C}$, inner semicontinuity of~$G$ in \ref{item:N5'} follows due to flows not being possible on~$D\backslash\widetilde{C}$. Further assuming~$x(0,1)\in\widetilde{C}\cup D$ and that~$x$ either flows or jumps after the first jump, inner semicontinuity of~$\widetilde{G}$ follows by a simple contradiction argument.
\end{pf}
}

\OnlyForAuto{See~\cite{arxiv} for a proof.} Roughly speaking, nominal inner well-posedness under~\ref{item:N1}-\ref{item:N4} can be guaranteed if every jump from a point~$\xi$ can be approximated either by another jump, or by a solution that flows for a short duration of time before jumping.

\begin{exmp}
\label{ex:great}
Consider the hybrid system~$\HS$ with data $=(C,F,D,G)$, where~$C=[0,1/2]$,~$D=[-1/2,1/2]$,~$F(x)=-1$ for all~$x\in\reals$, and~$G(x)=x$ if~$x\leq 0$ and~$G(x)=1$ otherwise. This system violates \ref{item:N5}-\ref{item:N6} due to discontinuity of~$G$, but satisfies~\ref{item:N1}-\ref{item:N4} and~\ref{item:N5'}-\ref{item:N6'}. Moreover, it can be verified to be nominally inner well-posed: nominal inner well-posedness away from the origin is obvious, and given a sequence~$\{\xi_i\}_{i=0}^{\infty}$ convergent to zero, the unique maximal solution from the origin can be approximated by letting,~$x_i(t,0)=\xi_i-t$ for all~$t\leq \xi_i$ and~$x_i(\xi,j)=0$ for all~$j\geq$ if~$\xi_i>0$, and~$x_i(0,j)=\xi_i$ otherwise.
\end{exmp}

%%%%%%%%%%%%%%%%%%%%%%%%%%%%%%%%%%%%%%%%%%%%%%%%%%%%%%%%%%%%%%%%%%%%%%%%%%%%%%%%%%%%%%%%%%

\subsection{Inner Well-Posed Perturbations}

As in nominal inner well-posedness, in developing sufficient conditions for inner well-posed perturbations,
terminal constraints arise naturally. In Definition~\ref{def:iwppertterm} below, while we make no assumptions on the family of terminal constraints~$\{X_{\delta}\}$, assuming that every point belonging to the set~$X$ is reached by a solution of~$\HS$ from~$S$, it is obviously necessary to have $X\subset\liminf_{\delta\to 0}X_{\delta}$ for~\ref{item:asterisk2} to hold. Note from the definition that a hybrid system~$\HS$ that is nominally inner well-posed with terminal constraint~$X$ can be viewed as an inner well-posed perturbation of itself with terminal constraints~$\{X_{\delta}=X\}$. \OnlyForAuto{See \cite[Example 23]{arxiv} for an example illustrating this concept.}

\begin{defn}
\label{def:iwppertterm}
A family of hybrid systems~$\{\HS_{\delta}\}$ with data~$\{(C_{\delta},F_{\delta},D_{\delta},G_{\delta})\}$ and terminal constraints~$\{X_{\delta}\}$ is an inner well-posed perturbation of~$\HS$ with terminal constraint~$X$ on a set~$S$ if~$S\cap(\closure(C)\cup D)\subset\liminf_{\delta\to 0}\closure(C_{\delta})\cup D_{\delta}$, and for every solution~$x$ of~$\HS$ originating from~$S$ and terminating on~$X$, the following holds.
\begin{enumerate}[label={($\diamond\diamond$)},leftmargin=*]
	\item \label{item:asterisk2} Given any sequence~$\{\delta_i\}_{i=0}^{\infty}\in(0,1)$ convergent to zero and any sequence~$\{\xi_{i}\}_{i=0}^{\infty}$ convergent to~$x(0,0)$ such that~$\xi_i\in\closure(C_{\delta_i})\cup D_{\delta_i}$ for all~$i\geq 0$, there exists~$\imath\geq 0$ such that the following holds: for every~$i\geq \imath$, there exists a solution~$x_i$ of~$\HS_{\delta_i}$ originating from~$\xi_i$ and terminating on~$X_{\delta_i}$ such that the sequence~$\{x_i\}_{i=0}^{\infty}$ is locally eventually bounded and graphically convergent to~$x$.
\end{enumerate} 
\end{defn}

\NotForAuto{
\begin{exmp}[Waypoint Navigation] Let~$(C_q,F_q)$ be a continuous-time system for all~$q\in\{0,1,\dots,N\}$ and~$\{p_q\}_{q=0}^{N+1}$ be a sequence of points such that for each~$q\leq N$, $p_q,p_{q+1}\in C_q$, and $(C_q,F_q)$ has a solution originating from~$p_q$ that terminates at~$p_{q+1}$. Consider the hybrid system~$\HS=(C,F,D,G)$ with state~$x=(z,q)$, where~$C = \cup_{q=0}^N \left(C_q\times\{q\}\right)$, $F(z,q)=F_q(z)\times\{0\}$ for all~$(z,q)\in C$, $D=\cup_{q=0}^N\{(p_{q+1},q)\}$, and $G(z,q)=(z,q+1)$ for all~$(z,q)\in\reals^n$. Note that~$\HS$ has a maximal solution~$x$ originating from~$(p_0,0)$ that can be interpreted as a waypoint navigation scenario with~$p_1,p_2,\dots,p_{N+1}$ the set of waypoints. If the continuous-time system~$(C_0,F_0)$ is not nominally inner well-posed at~$p_0$ with terminal constraint~$p_{1}$, clearly~$\HS$ is not nominally inner well-posed at~$(p_0,0)$ (with terminal constraint~$(p_{N+1},N+1)$). Hence, waypoint following is nonrobust; for some~$\tau\geq 0$ and~$\varepsilon>0$, there exists a sequence~$\{\xi_i\}_{i=0}\in C_0$ convergent to~$p_0$ such that for every~$i\geq 0$, no solution~$x_i$ of~$\HS$ from~$(\xi_i,0)$ is such that~$x$ and~$x_i$ are~$(\tau,\varepsilon)$-close~\cite[Proposition~3.18]{cdc2020}. Assuming uniqueness of solutions, this could imply, for example, that arbitrarily small perturbations to the initial conditions result in solutions not being able to reach the target~$(p_{N+1},N+1)$.

On the other hand, if each $(C_q,F_q)$ is nominally inner well-posed at~$p_q$ with terminal constraint~$S_{q+1}$ for some~$S_{q+1}\subset(C_q\cap C_{q+1})$ containing~$p_{q+1}$, one can construct the family of hybrid systems~$\{\HS_{\delta}\}$ with data~$\{(C_{\delta},F_{\delta},D_{\delta},G_{\delta})\}=\{(C,F,\widetilde{D},G)\}$, where the jump set~$\widetilde{D} := \cup_{q=0}^N (S_{q+1}\times\{q\})$. Then,~$\{\HS_{\delta}\}$ with terminal constraints~$\{X_{\delta}\}=\{S_{N+1}\times\{N+1\}\}$ is an inner well-posed perturbation of~$\HS$ with terminal constraint~$\{S_{N+1}\times\{N+1\}$ at~$(p_0,0)$. Note that for any~$q$, if~$(C_q,F_q)$ is nominally inner well-posed at~$p_q$, then it is also nominally inner well-posed at~$p_q$ with terminal constraint~$S_{q+1}$ for any relatively open subset of~$C_q\cap C_{q+1}$.
\end{exmp}
}

The following result generalizes Theorem~\ref{thrm:nomiwp}. The proof is in Appendix~\ref{sec:proof1}. Note that condition \eqref{eq:c5} and its analogue in \ref{item:N6p} can be viewed as an inner semi-continuity  property, with a restriction on the domain of convergence of the parameters.

\begin{thm}[Sufficient Conditions]
\label{thrm:nomiwppert}
Given a hybrid system~$\HS=(C,F,D,G)$, let~$\widetilde{C}$ be the set of all points where flows are possible. Similarly, given a family of hybrid systems~$\{\HS_{\delta}=(C_{\delta},F_{\delta},D_{\delta},G_{\delta})\}$, for every~$\delta\in(0,1)$, let~$\widetilde{C}_{\delta}$ be the set of all points where flows are possible for~$\HS_{\delta}$. Then,~$\{\HS_{\delta}\}$ is an inner well-posed perturbation of~$\HS$ if the following conditions hold:
\begin{enumerate}[label={(C\arabic*)},leftmargin=*]
	\item \label{item:N1p} for every~$x\in\widetilde{C}$, there exist~$r,\bar{\delta}>0$ such that for all~$\delta\in(0,\bar{\delta}]$, $(x+r\ball)\cap(\closure(C_{\delta})\cup D_{\delta})\subset\widetilde{C}_{\delta}$;
	\item \label{item:N2p} for every~$x\in D$, there exist~$r,\bar{\delta}>0$ such that for all~$\delta\in(0,\bar{\delta}]$, $(x+r\ball)\cap(\closure(C_{\delta})\cup D_{\delta})\subset\widetilde{C}_{\delta}\cup D_{\delta}$;
	\item \label{item:N3p} the family of continuous-time systems $\{(C_{\delta},F_{\delta})\}$ is an inner well-posed perturbation of the continuous-time system~$(C,F)$;
	\item \label{item:N4p} the family of continuous-time systems $\{(C_{\delta},F_{\delta})\}$ with terminal constraints~$\{D_{\delta}\}$ is an inner well-posed perturbation of the continuous-time system~$(C,F)$ with terminal constraint~$D$;
	\item \label{item:N5p} for every~$x\in D$,~$x\in\liminf_{\delta\to 0}D_{\delta}$, $x\notin\closure{C}$ implies $x\notin\limsup_{\delta\to 0}\closure C_{\delta}$, and
	\begin{equation}
	\label{eq:c5}
		G(x)\subset \liminf_{\substack{\delta\to 0,\, \xi \to x\\\xi \in D_{\delta}}}G_\delta(\xi);
	\end{equation}
	\item \label{item:N6p} for every~$x\in D$,
	\[
		G(x)\cap(\widetilde{C}\cup D)\subset \liminf_{\substack{\delta\to 0,\, \xi \to x\\\xi \in D_{\delta}}}G_\delta(\xi)\cap(\widetilde{C}_{\delta}\cup D_{\delta}).
	\]
\end{enumerate}
\end{thm}

Just like nominally inner well-posed hybrid systems, some of the sufficient conditions for inner well-posed perturbations turn out to be necessary. The proof resembles that of Theorem~\ref{prop:nec} and is omitted.

\begin{thm}[Necessary Conditions]
Given a hybrid system~$\HS=(C,F,D,G)$, let~$\widetilde{C}$ be the set of all points where flows are possible. Let~$\{\HS_{\delta}\}$ be an inner well-posed perturbation of~$\HS$ with data~$\{(C_{\delta},F_{\delta},D_{\delta},G_{\delta})\}$, and for every~$\delta\in(0,1)$, let~$\widetilde{C}_{\delta}$ be the set of all points where flows are possible for~$\HS_{\delta}$. Then,~\ref{item:N1p}-\ref{item:N2p} hold and~\ref{item:N3p}-\ref{item:N4p} hold if~$\closure C\subset\liminf_{\delta\to 0}\closure C_{\delta}$. Moreover 
\begin{enumerate}[label={(C\arabic*')},leftmargin=*]
\setcounter{enumi}{4}
\item for every~$x\in D$, if there exists a sequence~$\{\delta_i\}_{i=0}^{\infty}\in(0,1)$ convergent to zero and a sequence~$\{\xi_i\}_{i=0}^{\infty}$ convergent to~$x$ such that~$\xi_i\in D_{\delta_i}\backslash\widetilde{C}_{\delta_i}$, then
\begin{equation*}
		G(x)														\subset \liminf_{\substack{\delta\to 0,\, \xi \to x\\\xi \in D_{\delta}\backslash{\widetilde{C}_{\delta}}}}G_\delta(\xi),
\end{equation*}
\item for every~$x\in D$, if there exists a sequence~$\{\delta_i\}_{i=0}^{\infty}\in(0,1)$ convergent to zero and a sequence~$\{\xi_i\}_{i=0}^{\infty}$ convergent to~$x$ such that~$\xi_i\in D_{\delta_i}\backslash\widetilde{C}_{\delta_i}$, then
\begin{equation*}
		G(x)\cap(\widetilde{C}\cup D)		\subset \liminf_{\substack{\delta\to 0,\, \xi \to x\\\xi \in D_{\delta}\backslash{\widetilde{C}_{\delta}}}}G_\delta(\xi)\cap(\widetilde{C}_{\delta}\cup D_{\delta}).
\end{equation*}
\end{enumerate}
\end{thm}

We conclude the section with a proof of Theorem~\ref{thrm:nomiwp} using Theorem~\ref{thrm:nomiwppert}, as promised earlier. 

\begin{pf*}{Proof of Theorem~\ref{thrm:nomiwp}}
For each~$\delta\in(0,1)$, let~$C_{\delta}=C$,~$F_{\delta}=F$,~$D_{\delta}=D$, and~$G_{\delta}=G$. Then,~$\HS_{\delta}=\HS$ for all~$\delta\in(0,1)$. Similarly, let~$X_{\delta}=X$ for all~$\delta\in(0,1)$. Since \ref{item:N1p}-\ref{item:N6p} can be equivalently expressed as \ref{item:N1}-\ref{item:N6}, by Theorem~\ref{thrm:nomiwppert},~$\HS$ is nominally inner well-posed.
\end{pf*}

\section{Viability Conditions for Inner Well-Posedness}
\label{sec:viab}
To verify inner well-posedness of a hybrid system~$\HS$ with data~$(C,F,D,G)$ via Theorem~\ref{thrm:nomiwp}, one needs to be able to characterize continuous-time dynamics, in the sense that the set~$\widetilde{C}$ in~\eqref{eq:setflow} needs to be identified, and inner well-posedness of~$(C,F)$ (with terminal constraint~$D$) needs to be verified. This section addresses this by laying out viability conditions that characterize where flows are possible and ensure inner well-posedness for the continuous-time system~$(C,F)$. The conditions are then generalized so Theorem~\ref{thrm:nomiwppert} can be verify whether a family~$\{\HS_{\delta}\}$ is an inner well-posed perturbation of~$\HS$.

\subsection{Tangent Cones and Viability in Continuous Time}

In developing the viability conditions, we rely on two different types of tangent cones. Given a set~$S\subset\reals^n$ and a point~$x\in\reals^n$, the \textit{Bouligand tangent cone} (also known as the contingent cone) to~$S$ at~$x$ (\cite[Definition 5.12]{hybridbook} or \cite[Definition~6.1]{rockafellarwets}), denoted~$T_S(x)$, is the set of all~$v$ such that~$v=\lim_{i\to\infty}(x_i-x)/\tau_i$ for a sequence of points~$\{x_i\}_{i=0}^{\infty}\in S$ convergent to~$x$ and a positive sequence~$\{\tau_i\}_{i=0}^{\infty}$ convergent to zero. Similarly, the \textit{Dubovitsky-Miliutin tangent cone} to~$S$ at~$x$ \cite[Definition 4.3.1]{viability}, denoted~$M_S(x)$, is the set of all~$v$ for which there exist~$r,\bar{\delta}>0$ such that~$x+\delta w\in S$ for all~$\delta\in(0,\bar{\delta}]$ and~$w\in v+r\ball$. For every~$x\in\reals^n$,~$M_S(x)\subset T_S(x)$. For each~$x\in\partial S$,~$M_S(x)=\reals^n\backslash T_{\reals^n\backslash S}(x)$ \cite[Lemma~4.3.2]{viability}.

We frequently make use of the following lemmas. The first one is essentially~\cite[Lemma 5.26]{hybridbook}. The second one follows from \cite[Corollary 5.3.2]{viability} and the definition of Dubovitsky-Miliutin cone; see~\cite[Theorem~4.2~(c)]{cai}.

\begin{lem}[Viability with Outer Semicontinuity]
\label{lem:viaosc}
Given a continuous-time system~$(C,F)$, suppose that the flow set~$C$ is closed and~\ref{item:A2} holds. Let~$\widetilde{C}$ be the set of all points where flows are possible. Then,~$\interior C\subset\widetilde{C}$. Moreover, given~$x\in\partial C$, if $F(x)\cap T_C(x)$ is empty, then~$x\notin\widetilde{C}$, and if there exists~$r>0$ such that~$F(x')\cap T_C(x')$ is nonempty for all~$x'\in (x+r\ball)\cap \partial C$, then~$x\in\widetilde{C}$.
\end{lem}

\begin{lem}[Viability with Lipschitz Flow Maps]
\label{lem:viaosc2}
Let~$(C,F)$ be a continuous-time system and~$\widetilde{C}$ be the set of all points where flows are possible. Given~$x\in \closure C$, suppose that~$F(x)\cap M_{\interior C}$ is nonempty, and there exists an extension of~$F|_{\closure C}$ that is closed valued and Lipschitz on a neighborhood of~$x$. Then,~$x\in\widetilde{C}$. In particular, if the flow set~$C$ has a nonempty interior and the flow map~$F$ is locally Lipschitz and closed valued on~$\interior C$, then~$\interior C\subset\widetilde{C}$.
\end{lem}

%%%%%%%%%%%%%%%%%%%%%%%%%%%%%%%%%%%%%%%%%%%%%%%%%%%%%%%%%%%%%%%%%%%%%%%%%%%%%%%%%%%%%%%%%%

\subsection{Nominally Inner Well-Posed Systems}

For a continuous-time system $(C,F)$ with locally Lipschitz flow map~$F$, roughly speaking, nominal inner well-posedness is guaranteed if the data $(C,F)$ satisfies~\ref{item:A2} of the hybrid basic conditions, and the flow set~$C$ is ``locally contractive'' in the sense that if~$F(x)\cap T_C(x)$ is nonempty, then it must be that~$F(x')\subset M_{\interior C}(x')$ for all~$x'\in C$ close to~$x$. The regularity requirement, combined with the tangent cone condition,\footnote{Filippov's theorem assumes an open flow set. To invoke it for a solution~$x$ that has a point~$x(t)\in\partial C$, one needs to consider an extension of the flow map~$F$ to a neighborhood~$S$ of~$x(t)$. The tangent cone condition requires~$F$ to point to the interior of~$C$ so that extended solutions do not leave~$C$.} allows one to invoke Filippov's theorem~\cite[Theorem~5.3.1]{viability} to reach the~$(\tau,\varepsilon)$-closeness result assumed in \cite[Proposition~3.8]{cdc2020}, and conclude nominal inner-well posedness.

For nominal inner well-posedness with terminal constraint~$D$, additional conditions are required on the boundary of~$D$. One of these conditions is that~$F(x)$ should have a vector pointing towards~$\interior(C\cap D)$; see \ref{item:R3} and the second item of \ref{item:R4} in Theorem~\ref{thrm:viab} below. This ensures that any solution terminating at~$x$ can be extended to terminate at the interior of~$D$, and approximating sequences of solutions required for nominal inner well-posedness with terminal constraint~$D$ can be constructed using this extension. Alternatively, one can check that $D$ contains a relatively open subset of~$C$ that contains~$x$ (first item of \ref{item:R4}). If each vector of~$F(x)$ is pointing to the exterior of the flow set~$C$, then $D$ must contain a relatively open subset of~$\partial C$ that contains~$x$ (third item of \ref{item:R4}). Both of these guarantee that approximating solutions terminate on~$D$.

As in Section~\ref{sec:nomiwp}, the following theorem, which uses the conditions outlined above, is a special case of its upcoming counterpart for parametrized families. The proof is deferred to the end of the section.

\begin{thm}
\label{thrm:viab}
Given a continuous-time system~$(C,F)$, suppose that the flow set~$C$ is closed and~\ref{item:A2} holds. Moreover, suppose that the following hold.
\begin{enumerate}[label={(V\arabic*)},leftmargin=*]
	\item \label{item:R1}	For every~$x\in C$, there exists an extension of~$F|_{C}$ that is closed valued and Lipschitz on a neighborhood of~$x$.
	\item \label{item:R2}	For every~$x\in \partial C$ such that~$F(x)\cap T_C(x)$ is nonempty, there exists~$r>0$ such that $F(x')\subset M_{\interior C}(x')$ for all~$x'\in(x+r\ball)\cap\partial C$.
	\end{enumerate}
	Then, $(C,F)$ is nominally inner well-posed. If in addition, the following hold for a given set~$D$, then $(C,F)$ is nominally inner well-posed with terminal constraint~$D$.
	\begin{enumerate}[label={(V\arabic*)},leftmargin=*,resume]
	\item \label{item:R3} For every~$x\in D\cap(\interior (C)\cap\partial D)$,~$F(x)\cap M_{\interior D}(x)$ is nonempty.
	\item \label{item:R4} For every~$x\in D\cap(\partial C\cap \partial D)$, either of the following hold:
	\begin{itemize}[leftmargin=*]
		\item $(x+r\ball)\cap C\subset D$ for some~$r>0$;
		\item $F(x)\cap M_{\interior (C\cap D)}(x)$ is nonempty;
		\item $F(x)\cap T_C(x)$ is empty and $(x+r\ball)\cap \partial C\subset D$ for some~$r>0$.
	\end{itemize}
\end{enumerate}
\end{thm}

\NotForAuto{
\begin{figure}
\centering
\def\svgwidth{0.9\columnwidth}
	%% Creator: Inkscape inkscape 0.92.4, www.inkscape.org
%% PDF/EPS/PS + LaTeX output extension by Johan Engelen, 2010
%% Accompanies image file 'viability1n_3.pdf' (pdf, eps, ps)
%%
%% To include the image in your LaTeX document, write
%%   \input{<filename>.pdf_tex}
%%  instead of
%%   \includegraphics{<filename>.pdf}
%% To scale the image, write
%%   \def\svgwidth{<desired width>}
%%   \input{<filename>.pdf_tex}
%%  instead of
%%   \includegraphics[width=<desired width>]{<filename>.pdf}
%%
%% Images with a different path to the parent latex file can
%% be accessed with the `import' package (which may need to be
%% installed) using
%%   \usepackage{import}
%% in the preamble, and then including the image with
%%   \import{<path to file>}{<filename>.pdf_tex}
%% Alternatively, one can specify
%%   \graphicspath{{<path to file>/}}
%% 
%% For more information, please see info/svg-inkscape on CTAN:
%%   http://tug.ctan.org/tex-archive/info/svg-inkscape
%%
\begingroup%
  \makeatletter%
  \providecommand\color[2][]{%
    \errmessage{(Inkscape) Color is used for the text in Inkscape, but the package 'color.sty' is not loaded}%
    \renewcommand\color[2][]{}%
  }%
  \providecommand\transparent[1]{%
    \errmessage{(Inkscape) Transparency is used (non-zero) for the text in Inkscape, but the package 'transparent.sty' is not loaded}%
    \renewcommand\transparent[1]{}%
  }%
  \providecommand\rotatebox[2]{#2}%
  \newcommand*\fsize{\dimexpr\f@size pt\relax}%
  \newcommand*\lineheight[1]{\fontsize{\fsize}{#1\fsize}\selectfont}%
  \ifx\svgwidth\undefined%
    \setlength{\unitlength}{463.12569586bp}%
    \ifx\svgscale\undefined%
      \relax%
    \else%
      \setlength{\unitlength}{\unitlength * \real{\svgscale}}%
    \fi%
  \else%
    \setlength{\unitlength}{\svgwidth}%
  \fi%
  \global\let\svgwidth\undefined%
  \global\let\svgscale\undefined%
  \makeatother%
  \begin{picture}(1,0.55055953)%
    \lineheight{1}%
    \setlength\tabcolsep{0pt}%
    \put(0,0){\includegraphics[width=\unitlength,page=1]{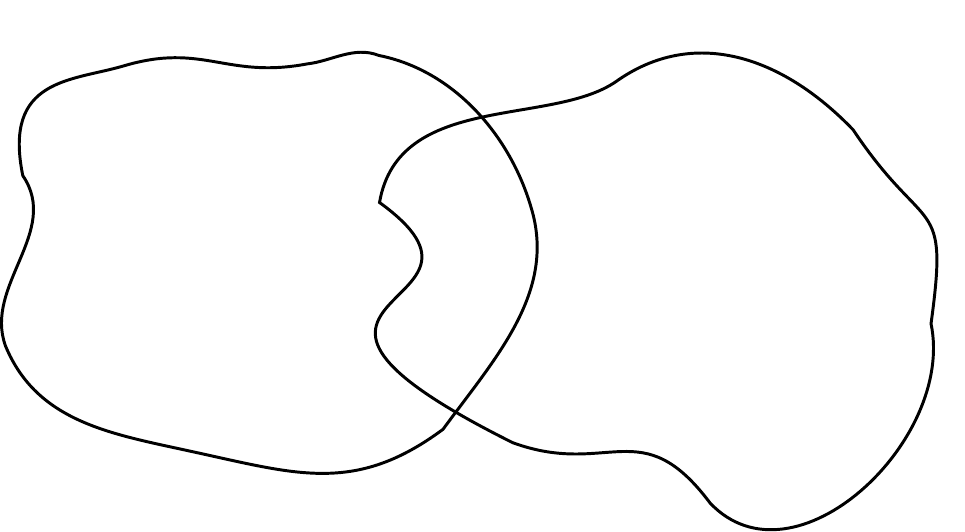}}%
    \put(0.87209028,0.48451715){\color[rgb]{0,0,0}\makebox(0,0)[lt]{\lineheight{1.25}\smash{\begin{tabular}[t]{l}$D$\end{tabular}}}}%
    \put(-0.00008877,0.50071143){\color[rgb]{0,0,0}\makebox(0,0)[lt]{\lineheight{1.25}\smash{\begin{tabular}[t]{l}$C$\end{tabular}}}}%
    \put(0,0){\includegraphics[width=\unitlength,page=2]{viability.pdf}}%
    \put(0.09060168,0.36277406){\color[rgb]{0,0,0}\makebox(0,0)[lt]{\lineheight{1.25}\smash{\begin{tabular}[t]{l}$F(x)$\end{tabular}}}}%
    \put(0.0457287,0.2812736){\color[rgb]{0,0,0}\makebox(0,0)[lt]{\lineheight{1.25}\smash{\begin{tabular}[t]{l}$T_C(x)$\end{tabular}}}}%
    \put(0,0){\includegraphics[width=\unitlength,page=3]{viability.pdf}}%
  \end{picture}%
\endgroup%

	\caption{Illustration of~\ref{item:R2}, given a point~$x\in\partial C$. $M_{\interior C}(x)$ could be interpreted as the set of vectors pointing to~$\interior C$ at~$x$, and~$T_C(x)$ as the set of vectors either pointing to~$\interior C$ or is tangential to~$\partial C$ at~$x$. Thus, for the point~$x$ in the figure,~$F(x)\cap T_C(x)$ is nonempty and~$F(x)\subset M_{\interior C}(x)$. When~$C$ is closed and~\ref{item:A2} holds, since~$F$ is upper semicontinuous relative to~$C$ and the tangent cone in the figure can be seen to vary smoothly around~$x$ (due to smoothness of~$\partial C$), it can be observed that $F(x')\subset M_{\interior C}(x')$ nearby~$x$.}
\end{figure}
}

The conditions above can also be used to verify \ref{item:N1}-\ref{item:N2} of Theorem~\ref{thrm:nomiwp}.
\begin{thm}
\label{thrm:viab2}
Given a continuous-time system~$(C,F)$, suppose that the flow set~$C$ is closed and~\ref{item:A2} and~\ref{item:R2} hold. Let~$\widetilde{C}$ be the set of all points from where flows are possible. Then,~$\widetilde{C}=\interior(C)\cup\{x\in\partial C: F(x)\cap T_C(x)\neq\varnothing \}$. Moreover,
\begin{itemize}[leftmargin=*]
	\item \ref{item:N1} holds if and if only if for every~$x\in\partial C$ such that~$F(x)\cap T_C(x)$ is nonempty, there exists~$r>0$ such that~$(x+r\ball)\cap D \subset C$;
	\item \ref{item:N2} holds if and if only if for every~$x\in D\cap(\partial C\cap \partial D)$ such that~$F(x)\cap T_C(x)$ is empty, there exists~$r>0$ such that $F(x')\cap T_C(x')$ is nonempty for all~$x'\in(x+r\ball)\cap ((\partial C)\backslash D)$. In particular, \ref{item:N2} holds if \ref{item:R4} holds.
\end{itemize}
\end{thm}
\begin{pf}
By Lemma~\ref{lem:viaosc}, the set $\widetilde{C}$ contains the interior of~$C$, as well as those points~$x$ on the boundary of~$C$ such that $F(x)\cap T_C(x)$ is nonempty, due to~\ref{item:R2} and the fact that~$M_{\interior C}(x')\subset T_C(x')$ for all~$x'$ on the boundary of~$C$. Note also that due to Lemma~\ref{lem:viaosc}, nonemptiness of~$F(x)\cap T_C(x)$ is necessary for any~$x\in\partial C$ to belong to~$\widetilde{C}$. This proves the equality regarding~$\widetilde{C}$, and the statements regarding~\ref{item:N1}-\ref{item:N2} follow immediately, since~$\widetilde{C}$ is relatively open in~$C$, as a result of \ref{item:R2}. The statement regarding~\ref{item:N2} and~\ref{item:R4} is also due to the Dubovitsky-Miliutin tangent cone being a subset of the Bouligand tangent cone on the boundary of~$C$.
\end{pf}

In~\cite[Example~1.2]{cdc2020}, we verify nominal inner well-posedness of a sample-and-hold control system using Theorems~\ref{thrm:nomiwp},~\ref{thrm:viab}, and~\ref{thrm:viab2}: the conditions of Theorem~1.1 therein used to prove nominal inner well-posedness are derived precisely by the conditions of Theorems~\ref{thrm:nomiwp},~\ref{thrm:viab}, and~\ref{thrm:viab2} here. Note that when the control law~$\kappa$ is continuous and the plant vector field~$f$ is Lipschitz (globally) for fixed~$u$, nominal inner well-posedness of sample-and-hold control can also be certified by Theorem~\ref{thrm:nomout} and \IfAuto{\cite[Proposition 8]{arxiv}}{Proposition~\ref{prop:ionwp}}, and the Picard-Lindel{\"o}f theorem.

%%%%%%%%%%%%%%%%%%%%%%%%%%%%%%%%%%%%%%%%%%%%%%%%%%%%%%%%%%%%%%%%%%%%%%%%%%%%%%%%%%%%%%%%%%

\subsection{Inner Well-Posed Perturbations}

To generalize Theorem~\ref{thrm:viab}, the following basic conditions are imposed on the parametrized family of continuous-time systems~$\{(C_{\delta},F_{\delta})\}$ with terminal constraints~$\{D_{\delta}\}$. Below, \ref{item:isccon} can be interpreted as the family of mappings $\{F_{\delta}\}$ ``uniformly converging to $F$'' over compact sets.
\begin{enumerate}[label={(P\arabic*)},leftmargin=*]
	\item \label{item:HBCeps} For every $\delta\in(0,1)$, the flow set $C_{\delta}$ is closed, the flow map~$F_{\delta}$ is locally bounded and outer semicontinuous relative to~$C_{\delta}$, and~$C_{\delta}\subset \dom F_{\delta}$. Furthermore, for every~$\delta\in(0,1)$ and~${x\in C_{\delta}}$, the set~${F_{\delta}(x)}$ is convex.
	\item \label{item:Cdelta} For every $\delta\in(0,1)$, $C_{\delta}\supset C$. 
	\item \label{item:isccon} For every compact set~$K\subset C$ and every~$\varepsilon>0$, there exists $\bar{\delta}>0$ such that $F(x)\subset F_{\delta}(x)+\varepsilon\ball$ for all~$x\in K$ and~$\delta\in(0,\bar{\delta}]$.
	\item \label{item:Ddelta} For every $\delta\in(0,1)$, $D_{\delta}\supset  C\cap \interior D$.
\end{enumerate}

The following theorem generalizes Theorem~\ref{thrm:viab} to establish parametrized families of continuous-time systems as inner well-posed perturbations of a given continuous-time system. The proof can be found in Appendix~\ref{sec:proof2}.

\begin{thm}
\label{thrm:viabpert}
Given a continuous-time system~$(C,F)$, suppose that the flow set~$C$ is closed and~\ref{item:A2} holds. Let~$\{(C_{\delta},F_{\delta})\}$ be a family of continuous-time systems satisfying~\ref{item:HBCeps}-\ref{item:isccon}, and suppose that the following hold.
\begin{enumerate}[label={(W\arabic*)},leftmargin=*]
	\item \label{item:Lip}For every~$x\in C$, there exist a neighborhood~$S$ of~$x$,~$\bar{\delta}>0$,~$L\geq 0$, and a family of mappings~$\{\widetilde{F}_{\delta}\}$ such that for every~$\delta\in(0,\bar{\delta}]$,~$\widetilde{F}_{\delta}$ is an extension of~$F_{\delta}|_{C_{\delta}}$ that is closed valued and Lipschitz on~$S$ with Lipschitz constant~$L$.
	 \item \label{item:DMdelta} For every compact~$K\subset \partial C$ such that~$F(x)\cap T_{C}(x)$ is nonempty for all~$x\in K$, there exist~$r,\bar{\delta}>0$ such that~$F_{\delta}(x')\subset M_{\interior C_{\delta}}(x')$ for all $x'\in(K+r\ball)\cap\partial C_{\delta}$ and $\delta\in(0,\bar{\delta}]$.
\end{enumerate}
Then, $\{(C_{\delta},F_{\delta})\}$ is an inner well-posed perturbation of~$(C,F)$. If, in addition, \ref{item:Ddelta} and the following hold for a given set~$D$ and a family of sets~$\{D_{\delta}\}$, then~$\{(C_{\delta},F_{\delta})\}$ with terminal constraints~$\{D_{\delta}\}$ is an inner well-posed perturbation of~$(C,F)$ with terminal constraint $D$.
\begin{enumerate}[label={(W\arabic*)},leftmargin=*,resume]
	\item \label{item:Dubovpert} For every~$x\in D\cap(\interior C\cap\partial D)$, either there exist~$r,\bar{\delta}>0$ such that $x+r\ball\subset D_{\delta}$ for all~$\delta\in(0,\bar{\delta}]$, or $F(x)\cap M_{\interior D}(x)$ is nonempty and~$F$ is Lipschitz on a neighborhood of~$x$.
	\item \label{item:last} For every~$x\in D\cap(\partial C\cap \partial D)$, there exist~$r,\bar{\delta}>0$, a neighborhood~$S$ of~$x$, and an extension~$\widetilde{F}$ of~$F$ such that either of the following hold.
		\begin{itemize}[leftmargin=*]
			\item $(x+r\ball)\cap C_{\delta}\subset D_{\delta}$ for all~$\delta\in(0,\bar{\delta}]$.
			\item $F(x)\cap M_{\interior(C\cap D)}(x)$ is nonempty and~$\widetilde{F}$ is Lipschitz and closed valued on~$S$.
			\item $F(x)\cap T_C(x)$ is empty,~$(x+r\ball)\cap (C_{\delta}\backslash \interior C)\subset D_{\delta}$ for all~$\delta\in(0,\bar{\delta}]$, either \ref{item:ay} or \ref{item:bee} below holds, and~\ref{item:Lip'} below holds.
			\begin{enumerate}[leftmargin=*]
				\item \label{item:ay} $\widetilde{F}$ is Lipschitz and closed valued on~$S$.
				\item \label{item:bee} $\widetilde{F}$ is locally bounded and outer semicontinuous relative to~$S$, $S\subset \dom \widetilde{F}$, and the set~$\widetilde{F}(x')$ is convex for all~$x'\in S$.
				\item \label{item:Lip'} There exists~$L\geq 0$ and a family of mappings~$\{\widetilde{F}_{\delta}\}$ such that for every~$\delta\in(0,\bar{\delta}]$,~$\widetilde{F}_{\delta}$ is an extension of~$F_{\delta}|_{C_{\delta}}$ that is closed valued and Lipschitz on~$S$ with Lipschitz constant~$L$. Moreover, for every compact~$K\subset S$ and every~$\varepsilon>0$, there exists $\bar{\delta}'>0$ such that $\widetilde{F}(x')\subset \widetilde{F}_{\delta}(x')+\varepsilon\ball$ for all~$x\in K$ and~$\delta\in(0,\bar{\delta}']$.
			\end{enumerate}
		\end{itemize}
\end{enumerate}
\end{thm}

\begin{rem}
Note that~\ref{item:Lip'} is essentially a combination of \ref{item:Lip} and \ref{item:isccon}. It allows the application of Filippov's theorem on solutions terminating on~$\partial C\cap\partial D$.
\end{rem}

\begin{exmp}[Nonlinear Oscillators]
\label{ex:osc}
Consider the harmonic oscillator described by the continuous-time system~$(C,F)$, where~$C$ is the unit circle on~$\reals^2$ and~$F(x)=(x_2,-x_1)$ for all~$x\in\reals^2$. By Theorem~\ref{thrm:nomout} and \IfAuto{\cite[Proposition 8]{arxiv}}{Proposition~\ref{prop:ionwp}}, this system is nominally inner well-posed. However, this cannot be verified using Theorem~\ref{thrm:viab}, since \ref{item:R2} is violated for every~$x\in \partial C=C$, as~$F(x)$ belongs to the tangent space to~$C$ at~$x$. Instead, we consider the family of continuous-time systems~$\{(C_{\delta},F_{\delta})\}$, where for every~$\delta\in(0,1)$, $C_{\delta}:=\{x: -\delta\leq x_1^2+x_2^2\leq 1+\delta\}$ and $F_{\delta}(x)=(x_2,-x_1)+(1-(x_1^2+x_2^2))(x_1,x_2)$ for all~$x\in\reals^2$. The perturbed flow map~$F_{\delta}$ is an extension of~$F|_C$. Clearly,~\ref{item:HBCeps}-\ref{item:isccon} hold. Outside~$C$,~$F_{\delta}(x)$ has the additional term $(1-(x_1^2+x_2^2))(x_1,x_2)$, which is orthogonal to~$F(x)$ and in the direction of~$C$, with its magnitude determined by that of~$x$. Therefore,~\ref{item:DMdelta} holds with arbitrary~$r,\bar{\delta}>0$. Moreover,~\ref{item:Lip} holds since~$F_{\delta}$, which is defined globally, does not depend on~$\delta$, and is locally Lipschitz. By Theorem~\ref{thrm:viabpert}, $\{C_{\delta},F_{\delta}\}$ is an inner well-posed perturbation of~$(C,F)$.
\end{exmp}

Here, we note that the simulation of the harmonic oscillator in Example~\ref{ex:osc} using the forward Euler method is not possible beyond a single integration step, due to the flow map being tangent to the unit circle. While we note that such problems can arise in the simulation of differential equations on manifolds, any solution of the harmonic oscillator can be approximated by solutions of the nonlinear oscillators~$\{(C_{\delta},F_{\delta})\}$ for small~$\delta$, which in turn can be simulated with desired precision.

Under some of the conditions of Theorem~\ref{thrm:viabpert}, \ref{item:N1p}-\ref{item:N2p} of Theorem~\ref{thrm:nomiwppert} simplify.

\begin{thm}
\label{thrm:viab2pert}
Given a continuous-time system~$(C,F)$, suppose that the flow set~$C$ is closed and~\ref{item:A2} holds. Let~$\widetilde{C}$ be the set of all points where flows are possible. Then, 
\[
	\interior(C)\subset\widetilde{C}\subset\interior(C)\cup\{x\in\partial C: F(x)\cap T_C(x)\neq\varnothing \}.
\]
 Similarly, given a family of continuous-time systems~$\{(C_{\delta},F_{\delta})\}$ satisfying~\ref{item:HBCeps} and~\ref{item:DMdelta}, for every~$\delta\in(0,1)$, let~$\widetilde{C}_{\delta}$ be the set of all points where flows are possible for~$\HS_{\delta}$. Then, $\interior C_{\delta} \subset \widetilde{C}_{\delta}$ for all~$\delta\in(0,1)$, and for every compact~$K\subset\partial C$ such that~$F(x)\cap T_{C}(x)$ is nonempty for all~$x\in K$, there exist~$r,\bar{\delta}>0$ such that $(K+r\ball)\cap\partial C_{\delta}\subset\widetilde{C}_{\delta}$ for all $\delta<\bar{\delta}$. Moreover, when~\ref{item:Cdelta} holds, the following hold:
\begin{itemize}[leftmargin=*]
	\item \ref{item:N1p} holds if for every~$x\in\partial C$ such that~$F(x)\cap T_C(x)$ is nonempty, there exists~$r,\bar{\delta}>0$ such that~$(x+r\ball)\cap D_{\delta} \subset C_{\delta}$ for all~$\delta\in(0,\bar{\delta}]$;
	\item assuming~$D_{\delta}\supset \interior D$ for all~$\delta \in(0,1)$,~\ref{item:N2p} holds if for every~$x\in D\cap(\partial C\cap \partial D)$ such that~$F(x)\cap T_C(x)$ is empty and every~$x\in D\cap(\partial D\backslash C)$, there exist~$r,\bar{\delta}>0$ such that~$(x+r\ball)\cap ((\partial C_{\delta})\backslash D_{\delta})\subset\widetilde{C}_{\delta}$. 
\end{itemize}
\end{thm}
\begin{pf}
The first two inclusions regarding $\widetilde{C}$ and $\widetilde{C}_{\delta}$ are due to Lemma~\ref{lem:viaosc}, the third inclusion is due to the Dubovitsky-Miliutin cone condition in \ref{item:DMdelta}, c.f. Theorem~\ref{thrm:viab2}. Since~$\interior C\subset\interior C_{\delta} \subset \widetilde{C}_{\delta}$ for all~$\delta\in(0,1)$ under~\ref{item:Cdelta}, for \ref{item:N1p}, it suffices to only check the boundary of~$C$. Then, given~$x\in\partial C$ with~$F(x)\cap T_C(x)$ nonempty, if there exist~$r,\bar{\delta}>0$ such that~$(x+r\ball)\cap D_{\delta} \subset C_{\delta}$ for all~$\delta\in(0,\bar{\delta}]$, the inclusion~$(x+r\ball)\cap D_{\delta} \subset C_{\delta}$ is enough to ensure~\ref{item:N1p}, since $(K+r\ball)\cap\partial (C_{\delta})\subset\widetilde{C}_{\delta}$ for all $\delta\in(0,\bar{\delta}]$ without loss of generality, due to the prior conclusion. The second item regarding~\ref{item:N1p} is then due to~$D_{\delta}\supset D$ (by assumption) and
~$\interior C\subset\interior C_{\delta} \subset \widetilde{C}_{\delta}$ under~\ref{item:Cdelta}.
\end{pf}

\begin{rem}
Under~\ref{item:Cdelta}, suppose that for every~$x\in (\partial D)\backslash C$, there exist~$r,\bar{\delta}>0$ such that~$(x+r\ball)\cap \partial C_{\delta}\subset D_{\delta}$ for all~$\delta\in(0,\bar{\delta}]$. In this case, to check that~\ref{item:N2p} holds, under additional regularity and tangent cone assumptions, one can rely on Lemmas~\ref{lem:viaosc} and~\ref{lem:viaosc2} to first fully characterize where flows are possible, which is also needed for~\ref{item:N6p}. Then, some of the conditions of Theorem~\ref{thrm:viabpert} can be used to guarantee~\ref{item:N2p}, as in Theorem~\ref{thrm:viab2}.
\end{rem}

\IfAuto
{
The section is concluded with the proof of Theorem~\ref{thrm:viab}. An example demonstrating the application of Theorems~\ref{thrm:nomiwppert}, \ref{thrm:viabpert}, and \ref{thrm:viab2pert} is in \cite[Example 35]{arxiv}.
}
{
The section is concluded with an example demonstrating the application of Theorems~\ref{thrm:nomiwppert}, \ref{thrm:viabpert}, and \ref{thrm:viab2pert}, along with the proof of Theorem~\ref{thrm:viab}. As mentioned before, we prove Theorem~\ref{thrm:viab} using Theorem~\ref{thrm:viabpert}.

\begin{exmp}[Bouncing Ball Revisited]
\label{ex:iwppertball}
Consider the bouncing ball system~$\HS=(C,F,D,G)$ in Example~\ref{ex:bouncingball}. This system is nominally inner well-posed, which can be verified by Theorem~\ref{prop:ionwp}. Consequently, the continuous-time dynamics of of the system is also nominally inner well-posed (with terminal constraint~$D$) by Proposition~\ref{prop:nec}, but this cannot be verified using Theorem~\ref{thrm:viab}, since \ref{item:R2} is violated when~$x$ is the origin. Instead, for some fixed~$r,c_1,c_2>0$, consider the family of hybrid systems~$\{\HS_{\delta}\}=\{(C_{\delta},F_{\delta},D_{\delta},G_{\delta}\}$, where for every~$\delta\in(0,1)$, $C_{\delta}=C\cup \{x: x_2\geq -c_2,c_1x_2\leq c_2x_1\}$, $D_{\delta}=\{x: \exists x'\in D,x\in(x'+r\ball)\cap C_{\delta}\} $, $G_{\delta}=G$, and
\[
	F_{\delta}(x)=
	\begin{cases}
		(0,-\gamma)								& c_1x_2> c_2x_1, x_2\leq 0\\
		(x_2-c_2x_1/c_1,-\gamma)	& c_1x_2\leq c_2x_1, x_1\leq 0\\
		(x_2,-\gamma)							& \text{otherwise}.
	\end{cases}
\]
Clearly, \ref{item:HBCeps}-\ref{item:Ddelta} hold, and since~$F_{\delta}$ does not depend on~$\delta$ and is locally Lipschitz, \ref{item:Lip} holds. For \ref{item:DMdelta}, one needs to only consider compact subsets of the nonnegative half of the~$x_2$-axis. Noting that~$F_{\delta}(x)\subset M_{\interior C_{\delta}}(x)$ for all $x'\in\partial (C_{\delta})$ and $\delta\in(0,1)$ such that~$x_2> -c_2$, it follows that \ref{item:DMdelta} holds. Moreover, \ref{item:Dubovpert} automatically holds since~$D\subset\partial C$, and by construction of~$D_{\delta}$, the first condition in \ref{item:last} holds for all~$x\in D$. Hence, Theorem~\ref{thrm:viabpert} applies. We also observe that Theorem~\ref{thrm:viab2pert} applies as well, by construction of $D_{\delta}\subset C_{\delta}$. Thus, \ref{item:N1p}-\ref{item:N4p} holds, and to be able to use Theorem~\ref{thrm:nomiwppert} to conclude inner well-posedness, it remains to prove \ref{item:N5p}-\ref{item:N6p}. The former holds trivially, as $G_{\delta}=G$ is a continuous function. Assume also that~$\lambda r < c_2$. Then, for every~$x=(x_1,x_2)\in D_{\delta}$, $G_{\delta}(x)=(0,-\lambda x_2)\in\widetilde{C}$ since~$-\lambda x_2>-c_2$ (due to~$F_{\delta}$ being locally Lipschitz and the inclusion~$F_{\delta}(x)\subset M_{\interior C_{\delta}}(x)$ for all $x'\in\partial (C_{\delta})$ with~$x_2> -c_2$, as a result of Theorem~\ref{thrm:viab2}). This observation and \ref{item:N5p} is sufficient to also conclude that \ref{item:N6p} holds. Therefore, Theorem~\ref{thrm:nomiwppert} applies;~$\{\HS_{\delta}\}$ is an inner well-posed perturbation of~$\HS$ if~$\lambda r < c_2$.
\end{exmp}

Here, it is noteworthy that in Example \ref{ex:iwppertball}, the Dubovitsky-Miliutin cone condition of \ref{item:last} can never hold since~$D$ has no interior. On the other hand, the tangent cone condition $F(x)\cap T_C(x)=\varnothing$ holds for every nonzero~$x\in D$, but it does not hold at the origin. Thus, for Theorem~\ref{thrm:viabpert} to apply, it is necessary for~$D_{\delta}$ to include the set~$r\ball\cap C_{\delta}$. In addition, unlike Example~\ref{ex:osc}, where the perturbed dynamics~$(C_{\delta},F_{\delta})$ tends to the nominal dynamics $(C,F)$ (as~$\delta$ approaches zero) and can be shown to be dominated by a~$\rho$-perturbation of the nominal system~$\HS$, the perturbed bouncing ball~$\HS_{\delta}$ does not tend to its nominal version~$\HS$. It is not obvious how this could be accomplished to satisfy~\ref{item:Lip} and~\ref{item:DMdelta} simultaneously. However, by letting~$r,c_1,c_2$ tend to zero, the nominal dynamics can be recovered.
}

\begin{pf*}{Proof of Theorem~\ref{thrm:viab}}
For each~$\delta\in(0,1)$, let~$C_{\delta}=C$,~$F_{\delta}=F$,~$D_{\delta}=D$, and~$G_{\delta}=G$. Then, \ref{item:HBCeps}-\ref{item:Ddelta} automatically hold when~$C$ is closed and \ref{item:A2} holds. Since\ref{item:Lip}-\ref{item:last} hold if \ref{item:R1}-\ref{item:R4} hold, by Theorem~\ref{thrm:viabpert}, $(C,F)$ is nominally inner well-posed (with terminal constraint $D$).
\end{pf*}

\section{Applications to Reachable Sets}
\label{sec:reachprop}
We apply the various well-posedness notions and results to reachable sets. In particular, we study the semicontinuous dependence of reachable sets on initial conditions, hybrid time, and perturbations by analyzing them within the framework of set-valued mappings. We primarily consider outer well-posed systems and inner well-posed perturbations. Corollaries for nominally outer/inner well-posed systems can be generated using~\cite[Theorem~4.1]{cdc2020}, and by treating nominally inner well-posed systems as an inner well-posed perturbations of themselves.

\begin{defn}[Reachable Set Mappings]
\label{def:reach}
Given a hybrid system~$\HS=(C,F,D,G)$, the reachable set mapping~$\reach_{\HS}:(\closure(C)\cup D)\times\realsgeq\times\nats\rightrightarrows\reals^n$ of~$\HS$ is the set-valued mapping associating with every~$x_0$,~$T$, and~$J$, the reachable set of~$\HS$ from~$x_0$ at time~$(T,J)$; $\reach_{\HS}(x_0,T,J):=\{x(T,J): x\in\sol_{\HS}(x_0), (T,J)\in\dom x\}$.
\end{defn}

An alternative formulation motivated by converse safety/invariance problems is found in~\cite{mohamed}, wherein the reachable set mappings collect the values of all solutions originating from~$x_0$ \textit{until} time~$(T,J)$. Yet another definition of reachable set mappings for hybrid inclusions can be found in~\cite[Section~6.3.2]{hybridbook}, which collects the values of all solutions originating from~$x_0$ for all~$(T,J)$ such that~$T+J\leq\tau$, given a time horizon parameter~$\tau\geq 0$.

\subsection{Semicontinuity of of Reachable Sets}

For a nominally outer well-posed pre-forward complete hybrid system~$\HS$, the reachable set mapping~$\reach_{\HS}$ is locally bounded and outer semicontinuous~\cite[Theorem~4.1]{cdc2020}. A consequence of this is that the reachable set of~$\HS$ from a compact set of initial conditions over a compact hybrid time horizon is compact~\cite[Proposition~4.2]{cdc2020}. Below, we generalize the former result to outer well-posed hybrid systems to account for perturbations. The proof follows the same steps, so it is not included. For a definition of~$\rho$-perturbations, recall Definition~\ref{def:rho}.

\begin{thm}[Outer Semicontinuity]
\label{thrm:roscoe}
Let~$\HS$ be a hybrid system. Given an initial condition~$x_0$, suppose that~$\HS$ is outer well-posed at~$x_0$ and pre-forward complete from~$x_0$. Then, given a continuous function~$\rho$, for every~$(T,J)\in\realsgeq\times\nats$, there exists~$\varepsilon>0$ such that the set~$\{\xi\in\reach_{\HS^{\delta\rho}}((x_0,T,J)+\varepsilon\ball): \delta \in(0,\varepsilon]\}=\reach_{\HS^{\varepsilon\rho}}((x_0,T,J)+\varepsilon\ball)$ is bounded, and
\begin{equation}
	\limsup_{\substack{\delta\to 0\\ (x_0',T',J') \to (x_0,T,J) }}\reach_{\HS^{\delta\rho}}(x_0',T',J') \subset\reach_{\HS}(x_0,T,J).
\label{eq:niceosc}
\end{equation}
\end{thm}

For a nominally inner well-posed system, inner semicontinuity of the reachable set mapping does not require pre-forward completeness. However, every point~$\xi$ belonging to the reachable set~$\reach_{\HS}(x_0,T,J)$ must be given by a maximal solution~$x$ that exhibits no jumps at ordinary time~$T$ and does not terminate at ordinary time~$T$; i.e., there exist~$x\in\sol_{\HS}(x_0)$ with~$x(T,J)=\xi$ and~$\varepsilon>0$ such that~$(T+\varepsilon,J)\in\dom x$ and~$(T,J-1)\notin\dom x$~\cite[Theorem~4.3]{cdc2020}. A similar conclusion holds for inner well-posed perturbations. Analogous conditions have been used for continuity of the reachability mapping in~\cite{mohamed}.

\begin{thm}[Inner Semicontinuity]
\label{thrm:neth}
Let~$\HS$ be a hybrid system. Given an initial condition~$x_0$ and~$(T,J)\in\realsgeq\times\nats$, suppose that for every~$\xi\in\reach_{\HS}(x_0,T,J)$, there exists~$x\in\sol_{\HS}(x_0)$ such that~$\xi=x(T,J)$ and~$T$ is not a jump time or the terminal ordinary time of~$x$. Let~$\{\HS_{\delta}\}$ be an inner well-posed perturbation of~$\HS$ at~$x_0$. Then,
\begin{equation}
	\reach_{\HS}(x_0,T,J)\subset\liminf_{\substack{\delta\to 0\\ (x_0',T',J') \to (x_0,T,J)}}\reach_{\HS_\delta}(x_0',T',J').
\label{eq:niceisc}
\end{equation}
\end{thm}
\begin{pf}
Pick~$x\in\sol_{\HS}(x_0)$ satisfying the required conditions. Take a nonempty interval~$[T_{\min},T_{\max}]$ with~$T<T_{\max}$ such that~$[T_{\min},T_{\max}]\times\{J\}\subset\dom x$, further assuming~$T>T_{\min}$ if~$T>0$, and noting that such an interval exists as~$T$ is not a jump time or the terminal ordinary time of~$x$. With no loss of generality, assume~$\dom x$ to be compact. Take any positive sequence~$\{\delta_i\}_{i=0}^{\infty}$ convergent to zero, any~$\{\xi'_i\}_{i=0}^{\infty}$ convergent to~$x_0$ with~$\xi'_i\in\closure(C_{\delta_i})\cup D_{\delta_i}$ for all~$i\geq 0$, and any~$\{T_i\}_{i=0}^{\infty}$ convergent to~$T$. Pick a sequence~$\{x'_i\}_{i=0}^{\infty}$ of hybrid arcs that is graphically convergent to~$x$ and locally eventually bounded, where, for each~$i\geq 0$,~$x'_i$ is a solution of~$\HS_{\delta_i}$ originating from~$\xi'_i$. By \IfAuto{\cite[Lemma 2]{arxiv}}{Lemma~\ref{lem:term}}, there exists~$\imath\geq 0$ such that~$(T_i,J)\in\dom x'_i$ for all~$i>\imath$. Due to local eventual boundedness of~$\{x'_i\}_{i=0}^{\infty}$, the sequence~$\{(T_i,J,x_i(T_i,J))\}_{i=\imath}^{\infty}$ is bounded, hence it has a convergent subsequence. By definition of graphical convergence, the limit of every such subsequence must belong to the graph of the outer graphical limit of~$\{x_i\}_{i=0}^{\infty}$, which is precisely~$x$, and since~$\lim_{i\to\infty}T_i=T$, the limit must be~$(T,J,x(T,J))$. It follows that~$\lim_{i\to\infty}(T_i,J,x_i(T_i,J))=(T,J,x(T,J))$. Therefore,~$x(T,J)\in\liminf_{i\to\infty}\reach_{\HS}(\xi_i,T_i,J)$.
\end{pf}

\begin{thm}[Continuity under Perturbations]
\label{thrm:ole}
Let~$\HS$ be a hybrid system, and given an initial condition~$x_0$, suppose that~$\HS$ is outer well-posed at~$x_0$ and pre-forward complete from~$x_0$. Let~$\{\HS_{\delta}\}$ be an inner well-posed perturbation of~$\HS$ at~$x_0$ that is dominated by a~$\rho$-perturbation of~$\HS$ for some continuous function~$\rho$. Moreover, given~$(T,J)\in\realsgeq\times\nats$, suppose that for every~$\xi\in\reach_{\HS}(x_0,T,J)$, there exists~$x\in\sol_{\HS}(x_0)$ such that~$\xi=x(T,J)$ and~$T$ is not a jump time or the terminal ordinary time of~$x$. Then,
\begin{equation}
	\reach_{\HS}(x_0,T,J)=\lim_{\substack{\delta\to 0\\ (x_0',T',J') \to (x_0,T,J)}}\reach_{\HS_\delta}(x_0',T',J'),
\label{eq:wee0}
\end{equation}
and the set~$\{\xi\in\reach_{\HS_{\delta}}((x_0,T,J)+\varepsilon\ball): \delta \in(0,\varepsilon]\}$ is bounded for some~$\varepsilon>0$.
\end{thm}

The proof of the continuity result above follows immediately by Theorems~\ref{thrm:roscoe} and~\ref{thrm:neth}. Counterexamples showing the necessity of the assumptions concerning the behavior of solutions in Theorem~\ref{thrm:neth} for inner semicontinuity are easy to construct: for the bouncing ball system Example~\ref{ex:bouncingball}, given a height parameter~$h>0$, consider~$x_0=(h,0)$ and~$T=\sqrt{2h/\gamma}$, for which~$\reach_{\HS}(x_0,T,0)=(0,\sqrt{2\gamma h})$. On the other hand, for every~$\varepsilon\in(0,h]$, $\reach_{\HS}(x_0-(\varepsilon,0),T,0)$) is empty since the first jump of the maximal solution from $x_0-(\varepsilon,0)$ occurs earlier than~$T$. 

\subsection{Continuous Approximations of Reachable Sets}

%%%%%%%%%%%%%%%%%%%%%%%%%%%%%%%%%%%%%%%%%%%%%%%%%%%%%%%%%%%%%

One of the main drawbacks of Theorem~\ref{thrm:neth} is that it requires partial knowledge of solutions. Moreover, if there exists~$\xi\in\reach_{\HS}(x_0,T,J)$ such that every maximal solution~$x$ from $x_0$ with $\xi=x(T,J)$ jumps or terminates at ordinary time $T$, at first glance, it is not clear whether a continuous, or even an inner semicontinuous approximation of~$\reach_{\HS}(x_0,T,J)$ is possible at all. 

With the following results, we show that to be able to continuously approximate the reachable set~$\reach_{\HS}(x_0,T,J)$, knowledge of solutions can be traded for class-$\K$ estimates characterizing lower semicontinuous dependence of solutions on initial conditions and perturbations (\IfAuto{\cite[Proposition 15]{arxiv}}{Proposition~\ref{prop:contout}}). Importantly, these results show that the reachable sets can be continuously approximated even when the solution-based conditions in Theorem~\ref{thrm:neth} are violated. The price to pay is that the approximation requires perturbed reachable sets over nontrivial continuous-time horizons, although as shown later, this can be relaxed.

\begin{thm}[Continuous-Time Inflations]
\label{thrm:cont1}
Let~$\HS$ be a hybrid system, and given a compact set~$K$, suppose that~$\HS$ is outer well-posed on~$K$ and pre-forward complete from~$K$. Let~$\{\HS_{\delta}=(C_{\delta},F_{\delta},D_{\delta},G_{\delta})\}$ be an inner well-posed perturbation of~$\HS$ on~$K$ that is dominated by a~$\rho$-perturbation of~$\HS$ for some continuous function~$\rho$. Then, given a compact set~$\T\subset\realsgeq\times\nats$, there exist class-$\K$ functions~$\alpha_1,\alpha_2$ such that for every~$x_0\in K$ and~$(T,J)\in \T$, the following holds: for every~$\varepsilon>0$ and~$\delta\in(0,\alpha_1(\varepsilon)]$, the set~$(x_0+\alpha_2(\varepsilon)\ball)\cap(\closure(C_{\delta})\cup D_{\delta})$ is nonempty if the reachable set~$\reach_{\HS}(x_0,T,J)$ is nonempty. Moreover,
\begin{multline}
	\lim_{\substack{\varepsilon,\delta\to 0,\, x_0' \to x_0\\ \varepsilon>0,\,\delta\in(0,\alpha_1(\varepsilon)]\\ x_0'\in x_0+\alpha_2(\varepsilon)\ball}}\reach_{\HS_{\delta}}(x_0',[\max\{0,T-\varepsilon\},T+\varepsilon],J)\\
	=\reach_{\HS}(x_0,T,J)
\label{eq:wee}
\end{multline}
\end{thm}
\begin{pf}
Pick~$\tau\geq 0$ such that~$T+J\leq\tau$ for all~$(T,J)\in\T$. By \IfAuto{\cite[Proposition 15, Remark 3.1]{arxiv}}{Proposition~\ref{prop:contout} and Remark~\ref{rem:notvac}}, there exist class-$\K$ functions~$\tilde{\alpha}_1,\alpha_2$ such that for all~$\varepsilon>0$, the following holds: for every~$x\in\sol_{\HS}(x_0)$ and $\delta\leq\tilde{\alpha}_1(\alpha_2(\varepsilon))$, the set~$S:=(x_0+\alpha_2(\varepsilon)\ball)\cap (\closure(C_{\delta})\cup D_{\delta})$ is nonempty, and for every~$x'_0\in S$, there exists a solution~$x'$ of~$\HS_{\delta}$ originating from~$x'_0$ such that~$x$ and~$x'$ are~$(\tau,\varepsilon)$-close. Letting~$\alpha_1:=\tilde{\alpha}_1\circ\alpha_2$, this is sufficient to prove the first statement. To prove~\eqref{eq:wee}, one needs to show that the right-hand side contains (respectively, is contained in) the outer (respectively, inner) limit of the left-hand side. That the right-hand side contains the outer limit of the left-hand side is an immediate result of Theorem~\ref{thrm:roscoe}. Now let~$x\in\sol_{\HS}(x_0)$ and without loss of generality, assume~$\dom x$ to be compact. Pick positive sequences~$\{\varepsilon_i\}_{i=0}^{\infty}$ and~$\{\delta_i\}_{i=0}^{\infty}$ tending to zero, where for every~$i\geq 0$,~$\delta_i\leq\alpha_1(\varepsilon_i)$. Take any sequence~$\{\xi'_i\}_{i=0}^{\infty}$ convergent to~$x_0$ such that~$\xi_i\in(x_0+\alpha_2(\varepsilon_i)\ball)\cap(\closure(C_{\delta_i})\cup D_{\delta_i})$ for all~$i\geq 0$. Let~$\{x'_i\}_{i=0}^{\infty}$ be a sequence of hybrid arcs that is graphically convergent to~$x$ and locally eventually bounded, where, for each~$i\geq 0$,~$x'_i$ is a solution of~$\HS_{\delta_i}$ originating from~$\xi'_i$. For large enough~$i$ (in the limit as~$i$ tends to infinity), there exists~$T_i$ such that~$|T_i-T|< \varepsilon_i$ and~$(T_i,J)\in\dom x'_i$. The rest of the proof is the same as that of Theorem~\ref{thrm:neth}.
\end{pf}
\begin{rem}
We emphasize nonemptiness of the set~$(x_0+\alpha_2(\varepsilon)\ball)\cap(\closure(C_{\delta})\cup D_{\delta})$, as it implies that the domain of the limit in~\eqref{eq:wee} is ``connected'', in the sense that given~$\varepsilon>0$ and~$\delta\leq\alpha_1(\varepsilon)$, there exists~$x_0'\in\closure(C_{\delta})\cup D_{\delta}$ that is $\alpha_2(\varepsilon)$-close to~$x_0$. Consequently, the perturbed reachable set~$\reach_{\HS_{\delta}}(x_0',[\max\{0,T-\varepsilon\},T+\varepsilon],J)$ is nonempty (for small~$\varepsilon>0$), provided the nominal reachable set~$\reach_{\HS}(x_0,T,J)$ is nonempty.
\end{rem}

The continuous approximation in Theorem~\ref{thrm:cont1} requires reachable set computations over hybrid time horizons corresponding to nontrivial continuous-time intervals. However, when flows are always possible after jumps, it suffices to consider computations over two hybrid times. The proof is similar to that of Theorem~\ref{thrm:cont1} and is omitted for brevity. In fact, the class-$\K$ functions~$\alpha_1,\alpha_2$ can simply be taken to be the ones in the proof of Theorem~\ref{thrm:cont1}. One needs to only note that due to the additional assumptions here, given~$\xi\in\reach_{\HS}(x_0,T,J)$, there exists~$x\in\sol_{\HS}(x_0)$ with~$(T,J)\in\dom x$ such that~$(T-\eta,J)\in\dom x$ or~$(T+\eta,J)\in\dom x$ for some~$\eta>0$.

\begin{thm}[Continuous-Time Doubling]
\label{thrm:cont2}
Under the conditions of Theorem~\ref{thrm:cont1}, if~$G(D)\subset\widetilde{C}$, where~$\widetilde{C}$ is the set of all points where flows are possible for~$\HS$, there exist class-$\K$ functions~$\alpha_1,\alpha_2$ such that for every~$x_0\in K$ and~$(T,J)\in \T$ satisfying~$T+J>0$ or~$x_0\in\widetilde{C}$,
\begin{multline}
	\reach_{\HS}(x_0,T,J)\\
	=\lim_{\substack{\varepsilon^-,\varepsilon^+,\delta^-,\delta^+\to 0,\,x_0^-,x_0^+ \to x_0 \\ 
	{\varepsilon},\delta^-,\delta^+>0,\,\delta\in(0,\alpha_1(\varepsilon)]\\ x_0^-\in x_0+\alpha_2(\varepsilon^-)\ball\\ x_0^+\in x_0+\alpha_2(\varepsilon^+)\ball}}S_{\delta^-}^-(x_0^-,\varepsilon^-)\cup S_{\delta^+}^+(x_0^+,\varepsilon^+),\\
		\begin{aligned}
		S_{\delta^-}^-(x_0^-,\varepsilon^-)&:=\reach_{\HS_{\delta^-}}(x_0^-,\max\{0,T-\varepsilon^-\},J),\\
		S_{\delta^+}^+(x_0^+,\varepsilon^+)&:=\reach_{\HS_{\delta^+}}(x_0^+,T+\varepsilon^+,J).
	\end{aligned}
\label{eq:wee2}
\end{multline}
where~${\varepsilon}:=\min\{\varepsilon^-,\varepsilon^+\}$ and~${\delta}:=\max\{\delta^-,\delta^+\}$.
\end{thm}

\subsection{Uniformity of Approximations}

It is worth stressing that uniformity of approximations can be achieved for all of the results derived in this section. For example, \eqref{eq:wee0} implies that for all~$\eta>0$ there exists~$\varepsilon>0$ such that for every~$\delta\in(0,\varepsilon]$, $x_0\in(\closure(C_{\delta})\cup D_{\delta})$, and~$(T',J')\in\realsgeq\times\nats$ satisfying~$|x'_0- x_0|\leq\varepsilon$ and~$|(T',J')-(T,J)|\leq\varepsilon$, $d(\reach_{\HS}(x_0,T,J),\reach_{\HS_{\delta}}(x'_0,T',J'))\leq\eta$, where~$d(.,.)$ is the Hausdorff distance. Similar relationships hold for \eqref{eq:niceosc}-\eqref{eq:wee2}, by~\cite[Proposition~5.12 and Exercise~5.13]{rockafellarwets}.

\section{Concluding Remarks}
\label{sec:conc}
Outer well-posedness of a hybrid system, defined originally in~\cite{hybridbook}, has paved the way for a fairly complete theory of hybrid systems as far as asymptotic stability is concerned. The notions of inner well-posedness introduced in this article, combined with outer well-posedness, ensure \textit{continuous} dependence of solutions initial conditions and perturbations and sets the stage for a more complete framework for hybrid systems:  we have already shown how these novel notions help lead to methods of continuously approximating reachable sets in this article, and proved consequent results (e.g continuous approximations of the value function) for optimal control problems for hybrid systems via reachability analysis (to appear, see~\cite{arxivReach} for an early version).

%%%%%%%%%%%%%%%%%%%%%%%%%%%%%%%%%%%%%%%%%%%%%%%%%%%%%%%%%%%%%

%\begin{ack}
%
%\end{ack}

%%%%%%%%%%%%%%%%%%%%%%%%%%%%%%%%%%%%%%%%%%%%%%%%%%%%%%%%%%%%%

\bibliographystyle{plain}
\bibliography{bibs/reach,bibs/arxiv}

%%%%%%%%%%%%%%%%%%%%%%%%%%%%%%%%%%%%%%%%%%%%%%%%%%%%%%%%%%%%%

\appendix

\IfAuto
{
\section{Proof of \texorpdfstring{Theorem \ref{thrm:nomiwppert}}{Theorem 17}}
}
{
\section{Proof of \texorpdfstring{Theorem \ref{thrm:nomiwppert}}{Theorem 24}}
}
\label{sec:proof1}
We begin by noting that by definition of inner well-posed perturbations,~\ref{item:N3p} implies~$\closure C\subset \liminf_{\delta\to 0} \closure C_{\delta}$. Hence, by \ref{item:N3p} and \ref{item:N5p}, $ \closure(C)\cup D\subset \liminf_{\delta\to 0} \closure (C_{\delta})\cup D_{\delta}$. In the rest of the proof, we show that the graphical convergence property in~\ref{item:asterisk} of Definition~\ref{def:iwppert} holds when the solution~$x$ therein has finitely many jumps. Then, using this property, we show that~\ref{item:asterisk} holds if~$x$ has infinitely many jumps, thanks to \IfAuto{\cite[Proposition 12 and 13]{arxiv}}{Propositions~\ref{prop:iscgraphical0pert} and~\ref{prop:graphicalisc0pert}}.

\subsection{The Case of Finite Number of Jumps}

Let~$x$ be a solution with~$J$ jumps, where~$J$ is finite. Let~$\{t_j\}_{j=1}^{J}$ be the jump times of~$x$ and~$t_0:=0$. For each~$j\in\{0,1,\dots,J\}$, let~$x^j$ be the restriction of~$x$ to the set of~$(s,i)\in\dom x$ with~$s+i\leq t_{j}+j$, and let~$x^{J+1}=x$. Since~$x$ has finitely many jumps, the proof will rely on an induction argument. That is, given a positive sequence~$\{\delta_i\}_{i=0}^{\infty}$ convergent to zero and a sequence~$\{\xi_{i}\}_{i=0}^{\infty}$ convergent to~$x(0,0)$, where $\xi_{i}\in\closure(C_i)\cup D_i$ for every~$i\geq 0$, we show that if there is a locally eventually bounded sequence of hybrid arcs~$\{x^j_{i}\}_{i=0}^{\infty}$ graphically convergent to~$x^j$ (with each~$x^j_i$ a solution of~$\HS_{\delta_i}$ originating from~$\xi_i$), then there is a sequence of hybrid arcs~$\{x^{j+1}_{i}\}_{i=0}^{\infty}$ graphically convergent to~$x^{j+1}$, where each~$x^{j+1}_i$ is a solution of~$\HS_{\delta_i}$ extending~$x^j_i$. Given~$j\leq J$, the induction hypothesis also assumes that there exists~$\imath\geq 0$ such that the following holds: i) for every~$i\geq \imath$, the interval~$I^j:=\{(t,j)\in\dom x^j_{i}\}$ is trivial (i.e., $x^j_{i}$ does not flow after jump~$j$); and ii) if~$j<J$ or if~$j=J$ and~$x^{J}\neq x^{J+1}$ (i.e., if~$x$ flows after jump~$j=J$),~$x^j_{i}$ terminates on~$\widetilde{C}_{\delta_i}\cup D_{\delta_i}$ for all~$i\geq\imath$.

The base case is obvious since~$x^0$ is trivial. In particular, one can simply take the sequence of trivial solutions~$\{x^0_i\}_{i=1}^{\infty}$, where each~$x^0_i$ originates from~$\xi_i$. Assuming that~$x$ is not trivial, note that~\ref{item:N1p}-\ref{item:N2p} ensure that these solutions terminate on~$\widetilde{C}_{\delta_i}\cup D_{\delta_i}$ for large~$i$. Now, given~$j\leq J$, let~$\{x^j_{i}\}_{i=0}^{\infty}$ be a locally eventually bounded sequence of hybrid arcs that is graphically convergent to~$x^j$, where~$x^j_i$ is a solution of~$\HS_{\delta_i}$ originating from~$\xi_i$ for all~$i\geq 0$. Note that there exists~$\imath\geq 0$ such that for every~$i\geq \imath$,~$\dom x^j_i$ is bounded (\IfAuto{\cite[Lemma 2]{arxiv}}{see Lemma~\ref{lem:term}}). Pass to this subsequence without relabeling. Without loss of generality, suppose that for each~$i\geq 0$,~$x^j_i$ has compact domain and let~$\xi'_i$ be the terminal point of~$x^j_i$. Recall that by \IfAuto{\cite[Lemma 2]{arxiv}}{the first bullet in Lemma~\ref{lem:term}}, $\lim_{i\to\infty}\xi'_i=x(t_j,j)$. To streamline the remainder of the proof, we rely on the following lemma, stated under the assumptions of Theorem \ref{thrm:nomiwppert}. It is used to extend~$x^j_{i}$ to~$x^{j+1}_{i}$ for each~$i\geq 0$. \IfAuto{Its proof can be found in \cite{arxiv}.}

\begin{lem}
\label{lem:forty}
Given a solution~$z$ of~$\HS$ with compact domain, suppose that~$z(T,1)\in\widetilde{C}\cup D$ and there exists no~$\varepsilon>0$ such that~$(T-\varepsilon,1)\in\dom z$, where~$(T,1)$ is the terminal time of~$z$. Let~$\{\delta_i\}_{i=0}^{\infty}$ be a positive sequence convergent to zero and~$\{\xi'_{i}\}_{i=0}^{\infty}$ be a sequence convergent to~$z(0,0)$ such that~$\xi'_{i}\in\closure(C_i)\cup D_i$ for all~$i\geq 0$. Then, there exist~$\imath\geq0$ and a locally eventually bounded sequence of hybrid arcs~$\{z_i\}_{i=\imath}^{\infty}$ graphically convergent to~$z$ such that the following holds for all~$i\geq\imath$: $z_i$ is a solution of~$\HS_{\delta_i}$ originating from~$\xi'_i$ and terminating on $\widetilde{C}_{\delta_i}\cup D_{\delta_i}$, and there exists no~$\varepsilon>0$ such that~$(T_i-\varepsilon,1)\in\dom z_i$, where~$(T_i,1)$ is the terminal time of~$z_i$.
\end{lem}
\NotForAuto{
\begin{pf}
Let~$z'$ be the truncation of~$z$ until hybrid time~$(T,0)$. We claim that there exist~$\imath\geq 0$ and a locally eventually bounded sequence of continuous arcs~$\{z'_i\}_{i=\imath}^{\infty}$ graphically convergent to~$z$ such that for every~$i\geq\imath$, $z'_i$ is a solution of~$\HS_{\delta_i}$ originating from~$\xi'_i$ and terminating on $D_{\delta_i}$. Indeed, if~$z'$ is nontrivial, then~$z(0,0)\in\widetilde{C}$, so existence of this sequence is guaranteed by~\ref{item:N1p} and~\ref{item:N4p}. If~$z'$ is trivial, then~$z(0,0)\in D$ and three cases are of interest. The first case is when~$\xi'_{i}\in D_{\delta_i}$ for large~$i$. This is straightforward to handle, as one can select the~$z'_i$'s to be trivial. The second case is when~$\xi'_{i}\in \closure C_{\delta_i}$ for large~$i$, which implies~$z(0,0)\in\limsup_{\delta\to 0}C_{\delta}$. Hence, by~\ref{item:N5p}, this implies that~$z(0,0) \in \closure (C)$, so existence of the sequence is guaranteed by~\ref{item:N4p}. The third case is the general case, in which $\{\xi'_{i}\}_{i=0}^{\infty}$ can be partitioned to two (infinite) subsequences, say~$\{\xi'_{i_k}\}_{k=0}^{\infty}$ and~$\{\xi'_{i_l}\}_{l=0}^{\infty}$, where~$\xi'_{i_k}\in \closure C_{\delta_{i_k}}$ for all~$k\geq 0$ and~$\xi'_{i_l}\in D_{\delta_{i_l}}$ for all~$l\geq 0$. Consequently, by combining the graphically convergent subsequences corresponding to~$\{\xi'_{i_k}\}_{k=0}^{\infty}$ and~$\{\xi'_{i_l}\}_{l=0}^{\infty}$, it is possible to construct~$\{z'_i\}_{i=\imath}^{\infty}$. 

Now, given the sequence $\{z'_i\}_{i=\imath}^{\infty}$, as a result of Lemma~\ref{lem:term}, we have that~$\lim_{i\to\infty}(T_i,z'_i(T_i,0))=z'(T,0)=z(T,0)$, where~$T_i$ is the terminal ordinary time of~$z'_i$ for all~$i\geq 0$. Hence, by~\ref{item:N6p}, there exists~$\imath\geq 0$ and a sequence~$\{\xi'_{i}\}_{i=\imath}^{\infty}$ convergent to~$z(T,1)$ such that~$\varsigma_{i}\in G(z'(T_i,0))\cap (\widetilde{C}_{\delta_i}\cup D_{\delta_i})$ for all~$i\geq\imath$. Letting, for each~$i\geq \imath$, $z_i(t,0)=z'_i(t,0)$ for all~$t\leq T_i$ and~$z_i(T_i,1)=\varsigma_i$ , the proof is complete.
\end{pf}
}

Lemma~\ref{lem:forty} allows us to move from induction step~$j$ to $j+1$ if~$j<J-1$. That is, we let~$z$ be the solution of~$\HS$ corresponding to~$x$ from~$(t_j,j)$ to~$(t_{j+1},j+1)$, and construct the sequence~$\{z_i\}_{i=\imath}^{\infty}$ in the lemma. Then, for~$i<\imath$, we let~$x^{j+1}_{i}=x^j_{i}$, otherwise we extend~$x^j_{i}$ to~$x^{j+1}_{i}$ by concatenating it with~$z_i$. It is then straightforward to show that the sequence~$\{x^{j+1}_{i}\}_{i=0}^{\infty}$ is locally eventually bounded and graphically convergent to~$x^{j+1}=x$. The same arguments apply for passing from step~$j=J-1$ to~$J$. The only delicate matter is the possibility that~$x^{J}= x^{J+1}$ (i.e.,~$J>0$ and~$x$ does not flow after jump~$J$, or~$J=0$ and~$x$ is trivial), as this implies that~$x(t_J,J)\notin\widetilde{C}$, and it could also be that~$x(t_J,J)\notin D$. In such a situation, the lemma above would need to be modified to account for~$z$ that terminates on~$\reals^n$ and to allow the~$z'_i$'s to also terminate on~$\reals^n$ using~\eqref{eq:c5} of~\ref{item:N5p} (instead of~\ref{item:N6p}).

Now, consider the case of~$j=J$. If~$x^{J}=x^{J+1}$ (in other words,~$x$ does not flow after jump~$J$), one can simply take the sequence~$\{x^{J+1}_{i}\}_{i=0}^{\infty}$ to be equal to~$\{x^J_{i}\}_{i=0}^{\infty}$. Otherwise, by the induction hypothesis, the terminal point~$\xi'_i$ belongs to~$\widetilde{C}_{\delta_i}\cup D_{\delta_i}$ for large~$i$. Moreover, since~$x(t_J,J)\in\widetilde{C}$, by~\ref{item:N1p} and \IfAuto{\cite[Lemma 2]{arxiv}}{Lemma~\ref{lem:term}}, there exists~$\imath\geq 0$ such that~$\xi'_i\in\widetilde{C}_{\delta_i}\subset\closure C_{\delta_i}$ for all~$i\geq \imath$. Similar to before, let~$z$ be the solution of~$\HS$ corresponding to~$x$ from~$(t_J,J)$ to~$(t_{J+1},J)$. By \ref{item:N3p}, there exists a locally eventually bounded sequence of continuous arcs~$\{z_{i}\}_{i=\imath}^{\infty}$ that is graphically convergent to~$z$, and for each~$i\geq \imath$, $z_i$ is a solution of~$(C_{\delta_i},F_{\delta_i})$. Thus, for each~$i\leq \imath$, we take~$x^{J+1}_{i}=x^J_i$, otherwise, we take~$x^{J+1}_{i}$ as the concatenation of~$x^J_i$ with~$z_i$. Again, it can be shown that the sequence~$\{x^{J+1}_{i}\}_{i=0}^{\infty}$ is locally eventually bounded and graphically convergent to~$x_{J+1}=x$.

\subsection{When the Solution has Infinitely Many Jumps}

We have shown that any solution~$x$ with finitely many jumps satisfies~\ref{item:asterisk}. The case of solutions with infinitely many jumps follows from \IfAuto{\cite[Proposition 12 and 13]{arxiv}}{Propositions~\ref{prop:iscgraphical0pert} and~\ref{prop:graphicalisc0pert}}. In particular, if~$x$ is a solution with infinitely many jumps, then for every~$\tau\geq 0$, the restriction of~$x$ to all~$(t,j)\in \dom x$ with~$t+j\leq\tau$, denoted~$x^{\tau}$, satisfies the following property: given a positive sequence~$\{\delta_i\}_{i=0}^{\infty}$ convergent to zero and a sequence~$\{\xi_{i}\}_{i=0}^{\infty}$ convergent to~$x(0,0)$ such that~$\xi_i\in\closure(C_{\delta_i})\cup D_{\delta_i}$ for all~$i\geq 0$, for every~$i\geq 0$, there exists a solution~$x_i$ of~$\HS_{\delta_i}$ originating from~$\xi_i$ such that the sequence~$\{x_i\}_{i=0}^{\infty}$ is locally eventually bounded and graphically convergent to~$x^{\tau}$. Then, by \IfAuto{\cite[Proposition 13]{arxiv}}{Proposition~\ref{prop:graphicalisc0pert}}, for every~$\varepsilon>0$ and~$\tau\geq 0$, there exist~$r,\bar{\delta}>0$ such that for any~$\delta\leq \bar{\delta}$ and~$x'_0\in (x(0,0)+r\ball)\cap (\closure(C_{\delta})\cup D_{\delta})$, there exists a solution~$x'$ of~$\HS_{\delta_i}$ originating from~$x'_0$ such that~$x^{\tau}$ and~$x'$ (and therefore,~$x$ and~$x'$) are~$(\tau,\varepsilon)$-close. By \IfAuto{\cite[Proposition 12]{arxiv}}{Proposition~\ref{prop:iscgraphical0pert}}, it follows that~$x$ satisfies~\ref{item:star}, completing the proof.

\IfAuto
{
\section{Proof of \texorpdfstring{Theorem \ref{thrm:viabpert}}{Theorem 23}}
}
{
\section{Proof of \texorpdfstring{Theorem \ref{thrm:viabpert}}{Theorem 30}}
}
\label{sec:proof2}
Observe that~$\closure C\subset \liminf_{\delta\to 0}\closure C_{\delta}$ by \ref{item:Cdelta}, and let~$x$ be a solution of~$(C,F)$ with closed graph. Given a positive sequence~$\{\delta_i\}_{i=0}$ convergent to zero, let~$\{\xi_{i}\}_{i=0}^{\infty}$ be a sequence convergent to~$x(0)$, where~$\xi_{i}\in C_{\delta_i}$ for all~$i\geq 0$. The solution has the property that if~$x(t)\in\partial C$ and~$F(x(t))\cap T_C(x(t))$ is empty for some~$t$, then~$t$ must be the terminal ordinary time of~$x$ due to Lemma~\ref{lem:viaosc}. Given~$\tau\geq 0$ such that the set~$\Omega:=\{x(t):t \leq \tau\}$ is bounded, let~$S$ be a neighborhood of~$\Omega$ and without loss of generality, by~\ref{item:Lip}, suppose that there exists~$L\geq 0$ such that for every~$\delta>0$, since~$\widetilde{F}_{\delta}$ is an extension of~$F_{\delta}$,~$F_{\delta}$ is closed valued and Lipschitz on~$S$ with Lipschitz constant~$L\geq 0$. Note that due to~\ref{item:HBCeps}, this property implies that~$F_{\delta}$ is upper semicontinuous on~$S$ and has nonempty and compact images on~$S$. By \ref{item:isccon} and Filippov's theorem (\cite[Theorem~5.3.1]{viability}), for every~$\varepsilon>0$, there exist~$r,\bar{\delta}>0$ such that for any~$\delta\in(0,\bar{\delta}]$ and every~$x'_0\in (x(0)+r\ball)\cap C_{\delta}$, there exists a solution~$x'$ of~$(S,F_{\delta})$ originating from~$x'_0$ such that~$x$ and~$x'$ are~$(\tau,\varepsilon)$-close---in fact, closeness can be quantified with the sup norm for the solutions \textit{and their derivatives}; namely $|x(t)-x'(t)|< \varepsilon$ for all~$t\leq\tau$ and $|\dot{x}(t)-\dot{x}'(t)|< \varepsilon$ for almost all~$t\leq\tau$.\footnote{Given small~$\delta>0$, straightforward application of Filippov's theorem would allow us to concur the existence of a solution~$x'$ such that~$x$ and~$x'$ are $(\tau,\varepsilon)$-close. The uniformity of the closeness result here over~$(0,\bar{\delta}]$ is due to the uniformity of the Lipschitz constant~$L$ in \ref{item:Lip} and the inclusion in \ref{item:isccon} over~$(0,\bar{\delta}]$; see (5.3) and the following equation in \cite{viability}.} If such a neighborhood does not exist, given a cover of~$\Omega$ coming from~\ref{item:Lip}, one can extract a finite subcover and apply Filippov's theorem repeatedly using the family of extended mappings $\widetilde{F}_{\delta}$ with~$\delta\in(0,\bar{\delta}]$.

\subsection{Inner Well-Posedness}

We consider two different cases to show that~$\{(C_{\delta},F_{\delta})\}$ is an inner well-posed perturbation of~$(C,F)$. We first show that the graphical convergence property required for inner well-posedness holds when the solution under question is bounded or complete, and has an additional tangent cone property. Afterward, we tackle the case where the solution escapes to infinity or fails to satisfy the aforementioned tangent cone condition.

\textbf{Case I:} Suppose that the solution~$x$ chosen above is bounded or complete, and has the property that there exists no~$t\geq 0$ such that~$x(t)\in\partial C$ and~$F(x(t))\cap T_C(x(t))$ is empty. If the range of~$x$ is contained in the interior of~$C$, then the $(\tau,\varepsilon)$-closeness property established above implies that for every~$i\geq 0$, there exists a solution~$x_i$ of~$(C_{\delta_i},F_{\delta_i})$ originating from~$\xi_i$ such that the sequence~$\{x_{i}\}_{i=0}^{\infty}$ is locally eventually bounded and graphically convergent to~$x$ by \IfAuto{\cite[Proposition 12]{arxiv}}{Proposition~\ref{prop:iscgraphical0pert}}, due to~\ref{item:Cdelta}. The same conclusion applies if the range of~$x$ intersects the boundary of~$C$: pick any~$\tau\geq 0$ and let~$K:=\{x(t)\in\partial C: t\leq \tau\}$. By \ref{item:DMdelta}, without loss of generality, there exists~$r>0$ such that~$F_{\delta}(\psi)\subset M_{\interior C_{\delta}}(\psi)$ holds for all~$\psi\in(K+r\ball)\cap\partial C_{\delta}$ and $\delta>0$. For every~$\delta>0$, let~$\eta_{\delta}:= \min_{t\in[0,\tau]}|x(t)|_{\partial (C_{\delta})\backslash (K+r\ball)}>0$. Then, there exists~$\eta$ such that~$\eta_{\delta}\geq \eta$ for all~$\delta> 0$, due to~\ref{item:Cdelta}. Take any~$\varepsilon\in(0,\min\{\eta/2,r\})$, any~$\delta'>0$, and any solution~$x'$ of~$(S,F_{\delta'})$ originating from~$C_{\delta'}$ such that~$x$ and~$x'$ are $(\tau,\varepsilon)$-close. Then,~$x'$ is a solution of~$(C_{\delta'},F_{\delta'})$. To show this, suppose the opposite, which implies that there exists~$t'\geq 0$ such that~$x'(t')\in \partial C_{\delta'}$, and a sequence~$\{t'_i\}_{i=0}^{\infty}$ convergent to~$t'$ such that~$t'_i>t'$ and~$x(t'_i)\notin C_{\delta'}$. Since~$x$ and~$x'$ are $(\tau,\varepsilon)$-close and~$\varepsilon< \min\{\eta/2,r\}$,~$x'(t')\in K+r\ball$. However, by~\cite[Theorem 4.3.6]{viability} and as a result of the Dubovitsky-Miliutin tangent cone condition, every solution~$\tilde{x}'$ of~$(S,F_{\delta'})$ from~$x'(t')$ stays in~$C_{\delta'}$ for some nonzero amount of time,\footnote{Application of~\cite[Theorem  4.3.6]{viability} is possible, thanks to the fact that Krasovskii regularization (\cite[Definition~4.13]{hybridbook}) of the mapping~$F_{\delta}$ is upper semicontinuous and has nonempty, compact, and convex values; see~\cite[Lemmas~5.15-5.16]{hybridbook}.} contradicting the existence of the sequence~$\{t'_i\}_{i=0}^{\infty}$. Since~$x'$ is a solution of~$(C_{\delta'},F_{\delta'})$, again, \IfAuto{\cite[Proposition 12]{arxiv}}{Proposition~\ref{prop:iscgraphical0pert}} is applicable.

\textbf{Case II:} Now recall that if there exists~$t\geq 0$ such that~$x(t)\in\partial C$ and~$F(x(t))\cap T_C(x(t))$ is empty,~$t$ is the terminal ordinary time of the solution~$x$ chosen above. Suppose that either~$x(T)\in\partial C$ and~$F(x(T))\cap T_C(x(T))$ is empty, where~$T$ is the terminal ordinary time of~$x$, or~$x$ has finite escape time~$T$. Let~$T_k:=Tk/(k+1)$ for every~$k\geq 0$ and fix~$\varepsilon>0$. We build the graphically convergent sequence with a recursive procedure:
\begin{itemize}[leftmargin=*]
	\item For~$k=0$, we take a locally eventually bounded sequence of solutions~$\{\tilde{x}_i^0\}_{i=0}^{\infty}$ (with each~$\tilde{x}_i^0$ a solution of~$(C_{\delta_i},F_{\delta_i})$ originating from~$\xi_i$) that is graphically convergent to the solution~$y^0$, where~$y^0$ corresponds to~$x$ from~$T_0=0$ to~$T_1$. Then, for each~$i\geq 0$, we let~$\hat{x}_i^0:=\tilde{x}_i^0$. Also, we let~$i_0=0$.
	\item For~$k\geq 1$, we take a locally eventually bounded sequence of solutions~$\{\tilde{x}_i^k\}_{i=0}^{\infty}$ (with each~$\tilde{x}_i^k$ a solution of~$(C_{\delta_i},F_{\delta_i})$ originating from the terminal point of~$\hat{x}_i^{k-1}$) that is graphically convergent to the solution~$y^k$, where~$y^k$ corresponds to~$x$ from~$T_k$ to~$T_{k+1}$. Then, for each~$i\geq k$, we let~$\hat{x}_i^k$ be the extension of~$\hat{x}_i^{k-1}$ with~$\tilde{x}_i^k$. Also, noting that~$\{\hat{x}_i^k\}_{i=0}^{\infty}$ is locally eventually bounded and converges to the truncation of $x$ up to~$T_{k+1}$, by~\cite[Theorem~5.25]{hybridbook}, we pick some~$i_k>i_{k-1}$ such that~$\hat{x}_{i}^k$ and the truncation of $x$ up to~$T_{k+1}$ are~$(T_{k+1},\varepsilon/k)$-close for all~$i\geq i_k$.
\end{itemize}
Lastly, for each~$i\geq 0$, we let~$x_i:=\hat{x}_{i}^k$, where~$k\geq 0$ is such that~$i_k\leq i <i_{k+1}$.

The limit of the sequence~$\{\dom {x}_i\}_{i=0}^{\infty}$ is precisely~$[0,T]$. Moreover, the outer and inner graphical limits of the sequence $\{ {x}_i\}_{i=0}^{\infty}$ are mappings~$M_{\text{out}}$ and~$M_{\text{in}}$, respectively, with the following properties: the domains of both mappings contain~$[0,T)$ and are subsets of~$[0,T]$, and the restriction of both mappings to~$[0,T)$ is equal to the restriction of~$x$ to~$[0,T)$. This can be observed using \IfAuto{\cite[Lemma 2]{arxiv}}{Lemma~\ref{lem:term}} by considering truncations of the~$x_i$'s; that is, given~$t< T$, for each~$i\geq 0$, we truncate~$x_i$ until some ordinary time~$t_i$ such that~$\{t_i\}_{i=0}^{\infty}$ tends to~$t$ and analyze the resulting sequence of solutions as follows. 
\begin{itemize}[leftmargin=*]
	\item If the solution~$x$ is bounded, since the inner and outer set limits are always closed, it follows that graphs of~$M_{\text{out}}$ and~$M_{\text{in}}$ contain~$(T,x(T))$. In addition, by construction (in particular, due to $(T_k+\varepsilon/k,\varepsilon/k)$-closeness),~$\{x_i\}_{i=0}^{\infty}$ is locally eventually bounded and~$M_{\text{out}}$ (and therefore~$M_{\text{in}}$) is single valued at~$T$. Hence, the graphical limit is precisely~$x$.
	\item On the other hand, if~$x$ has finite escape time, $\{ x_i\}_{i=0}^{\infty}$ is not locally bounded, as the converse contradicts the fact that~$x$ is unbounded. To show that the graphical limit leads to a set-valued mapping~$M$ that is an extension of~$x$, note that by construction, $M_{\text{out}}(T)=M_{\text{in}}(T)=\varnothing$: assume $M_{\text{out}}(T)$ nonempty, and without loss of generality, take a sequence~$\{t_i\}_{i=0}^{\infty}$ convergent to~$T$ such that~$\{x_i(t_i)\}_{i=0}^{\infty}$ is convergent. By construction, there exists a sequence~$\{s_i\}_{i=0}^{\infty}$ such that~$\lim_{i\to\infty} |t_i-s_i|=\lim_{i\to\infty} |x_i(t_i)-x(s_i)|=0$. Consequently,~$\{x(s_i)\}_{i=0}^{\infty}$ is convergent, which is a contradiction, as~$x$ has finite escape time~$T$. This proves that $\{(C_{\delta},F_{\delta}\})$ is an inner well-posed perturbation of~$(C,F)$, i.e.,~\ref{item:N3p} holds.
\end{itemize}

\subsection{Inner Well-Posedness with Terminal Constraints}

Let~$x$ be a solution of $(C,F)$ terminating on~$D$, with terminal ordinary time~$T$. Three cases are of interest.

\textbf{Case I:} If~$x(T)\in \interior D$, then the graphical convergence property for~$x$ holds without any additional conditions by~\ref{item:N3p} and \IfAuto{\cite[Lemma 2]{arxiv}}{Lemma~\ref{lem:term}}, as a result of \ref{item:Ddelta}.

\textbf{Case II:} Suppose $x(T)\in\interior (C)\cap\partial (D)$. If there exist~$r,\bar{\delta}>0$ such that $x(T)+r\ball\subset D_{\delta} $ for all~$\delta\leq\bar{\delta}$ (see~\ref{item:Dubovpert}), the graphical convergence property for~$x$ holds without any additional conditions by~\ref{item:N3p} and \IfAuto{\cite[Lemma 2]{arxiv}}{Lemma~\ref{lem:term}}. Otherwise, if~$F(x(T))\cap M_{\interior D}(x(T))$ is nonempty and~$F$ is closed valued and Lipschitz on a neighborhood of~$x(T)$ by~\ref{item:Dubovpert}, we extend~$x$ to another solution, say~$x'$, such that~$x'(t)\in \interior(C)\cap\interior (D)$ for all~$t\in(T,T']$---this is possible by Lemma~\ref{lem:viaosc2}, where~$T'>T$ is the terminal ordinary time of~$x'$. For every~$i\geq 0$, let~$\tilde{x}_i$ be a solution of~$(C_{\delta_i},F_{\delta_i})$ originating from~$\xi_i$ such that that the sequence~$\{\tilde{x}_{i}\}_{i=0}^{\infty}$ is locally eventually bounded and graphically convergent to~$x'$. Note that there exists~$\imath\geq 0$ such that for every~$i\geq\imath$, $\tilde{x}_i$ terminates on~$(\interior (C) \cap\interior (D))\subset D_{\delta_i}$ at terminal ordinary time~$\widetilde{T}_i>T$, by~\ref{item:Ddelta}. For each~$i\geq\imath$, pick~$T_i\in [T,\widetilde{T}_i]$ with~$\tilde{x}_i(T_i)\in D_{\delta_i}$ such that the sequence~$\{T_i\}_{i=\imath}^{\infty}$ tends to~$T$, and let~$x_i$ be the truncation of~$\tilde{x}_i$ until ordinary time~$T_i$. Existence of such~$\{T_i\}_{i=\imath}^{\infty}$ is due to the set limit being closed: local boundedness and graphical convergence of~$\{\tilde{x}_{i}\}_{i=0}^{\infty}$ implies that~$\lim_{i\to\infty}\mathcal{T}_i\supset [T,T']$, where~$\mathcal{ T}_i:=\{t\in[T,\widetilde{T}_i]:\tilde{x}_i(t)\in D_{\delta_i}\}$, due to~\ref{item:Ddelta} and the fact that~$x'(t)\in \interior D$ for all~$t\in(T,T']$. Hence~$\{x_{i}\}_{i=\imath}^{\infty}$ is graphically convergent to~$x$.

\textbf{Case III:} Finally, suppose that~$x(T)\in\partial (C)\cap\partial (D)$. If there exist~$r,\bar{\delta}>0$ such that~$(x(T)+r\ball)\cap C_{\delta}\subset D_{\delta}$ for all~$\delta\leq\bar{\delta}$, we take any locally eventually bounded sequence~$\{x_i\}_{i=0}^{\infty}$ graphically convergent to~$x$, where each~$x_i$ is a solution of~$(C_{\delta_i},F_{\delta_i})$ originating from~$\xi_i$. Then,~$x_i$ terminates on~$D_{\delta_i}$ for large enough~$i$. If~$F(x(T))\cap M_{\interior(C)\cap\interior(D)}(x(T))$ is nonempty, then~$x$ can be extended to a solution~$x'$ such that~$x'(t)\in\interior(C)\cap\interior(D)$ for all~$t\geq T$, as in Case II. The fact that~$x'$ stays in~$\interior(C)\cap\interior (D)$ is due to the Dubovitsky-Miliutin tangent cone condition in~\ref{item:last} and the invariance property outlined in~\cite[Theorem 4.3.6]{viability}. The locally eventually bounded graphically convergent sequence is then constructed as in Case II, in the same manner. 

Otherwise, if~$F(x(T))\cap T_C(x(T))$ is empty, due to either of the regularity properties of~$\widetilde{F}$ in~\ref{item:last},~$x$ can be extended as in Case II to another function with compact domain. However, this time, the extension~$x'$ is not a solution of~$(C,F)$, and there exists a positive sequence~$\{s_i\}_{i=0}^{\infty}$ convergent to zero such that~$x'(T+s_i)\notin C$ for all~$i\geq 0$, due to the tangent cone condition~$F(x(T))\cap T_C(x(T))=\varnothing$ (this precludes existence of solutions from~$x(T)$ by Lemma~\ref{lem:viaosc}). We use similar arguments as before. We first note that due to \ref{item:Lip'}, Filippov's theorem (\cite[Theorem~5.3.1]{viability}) is again applicable, this time for the extended solution~$x'$. This implies that there exists a locally eventually bounded sequence of continuous arcs~$\{\tilde{x}_{i}\}_{i=0}^{\infty}$ graphically convergent to~$x'$, where for each~$i\geq 0$, $\tilde{x}_{i}$ is an extension of a solution~$x'_i$ of~$(C_{\delta_i},F_{\delta_i})$ originating from~$\xi_i$, and the sequence~$\{x'_{i}\}_{i=0}^{\infty}$ converges to~$x$. For each~$i\geq 0$, let~$T'_i$ be the terminal time of~$x'_i$, and pick~$T_i\geq T'_i$ such that the following hold: 1)~$\tilde{x}_i(t)\in C_{\delta_i}$ for all~$t\leq T_i$, 2)~$\tilde{x}_i(T_i)\in C_{\delta_i}\backslash \interior C$, and 3) the sequence~$\{T_i\}_{i=0}^{\infty}$ converges to~$T$. Given~$i\geq 0$, let~$x_i$ be the truncation of~$\tilde{x}_i$ until ordinary time~$T_i$. As in the case $x(T)\in\interior (C)\cap\partial (D)$, the sequence~$\{x_{i}\}_{i=0}^{\infty}$ is graphically convergent to~$x$ and by \IfAuto{\cite[Lemma 2]{arxiv}}{Lemma~\ref{lem:term}}, there exists~$\imath\geq 0$ such that~$\{x_i(T_i,0)\}_{i=\imath}^{\infty}$ tends to~$x(T,0)$ and $x_i(T_i,0)\in D_{\delta_i}$ for all~$i\geq \imath$.

\NotForAuto
{
\section{Outer Well-Posedness Formalized}
\label{sec:defs}
\begin{defn}{\cite[Definition 3.2]{cdc2020}}
\label{def:nowp}
A hybrid system~$\HS$ is said to be \textit{nominally outer well-posed} on a set~$S\subset\reals^n$ if for every graphically convergent sequence of solutions~$\{x_i\}_{i=0}^{\infty}$ of~$\HS$ satisfying~$\lim_{i\to\infty}x_i(0,0)=:x_0\in S$, the following holds:
	\begin{enumerate}[label={(\alph*)},leftmargin=*]
		\item \label{item:heya} if the sequence~$\{x_i\}_{i=0}^{\infty}$ is locally eventually bounded, then the graphical limit~$x$ is a solution of~$\HS$ originating from~$x_0$;
		\item \label{item:heyb} if the sequence~$\{x_i\}_{i=0}^{\infty}$ is not locally eventually bounded, then there exists~$(T,J)\in\realsgeq\times\nats$ such that~$x=M|_{\dom M\cap([0,T)\times\{0,1,\dots,J\})}$ is a solution of~$\HS$ originating from~$x_0$ that escapes to infinity at time~$(T,J)$, where~$M$ is the graphical limit of~$\{x_i\}_{i=0}^{\infty}$.
	\end{enumerate}
\end{defn}

\begin{defn}[Outer Well-Posedness]
\label{def:owp}
A hybrid system~$\HS$ is said to be \textit{outer well-posed} on a set~$S\subset\reals^n$ if for every continuous function~$\rho$, every positive sequence~$\{\delta_i\}_{i=0}^{\infty}$ convergent to zero, and every graphically convergent sequence of hybrid arcs~$\{x_i\}_{i=0}^{\infty}$ such that~$x_i$ is a solution of~$\HS^{\delta_i\rho}$ and~$\lim_{i\to\infty}x_i(0,0)=:x_0\in S$, \ref{item:heya}-\ref{item:heyb} in Definition~\ref{def:nowp} hold.
\end{defn}
}

\end{document}